\documentclass[10pt]{amsart}
\usepackage{latexsym,enumerate}
\usepackage{color}
\usepackage[latin1]{inputenc}
\usepackage{amsmath,amsthm,amsopn,amstext,amscd,amsfonts,amssymb}
\usepackage[dvips]{graphicx}
\usepackage[active]{srcltx}
\newcommand{\R}{{I\!\!R}}
\newcommand{\N}{{I\!\!N}}

\newcommand{\car}{{\raise2pt\hbox{$\chi$}}}

\newcommand{\p}{\partial}

\newcommand{\Om}{\Omega}

\newcommand{\re}{{I\!\!R}}
\newcommand{\ren}{\re^N}

\newcommand{\dyle}{\displaystyle}

\newcommand{\io}{\int\limits_{\O}}

\newcommand{\sob}{W_0^{1,2}}

\newcommand{\Di}{\text{div}}

\renewcommand{\a }{\alpha }

\newcommand{\inn}{\mbox{ in }}
\newcommand{\onn}{\mbox{ on }}
\renewcommand{\b }{\beta }
\renewcommand{\d }{\delta }
\newcommand{\D }{\Delta }
\newcommand{\e }{\varepsilon }
\newcommand{\eps}{\varepsilon }

\newcommand{\g }{\gamma}

\renewcommand{\l }{\lambda }

\newcommand{\n }{\nabla }

\newcommand{\s }{\sigma }

\renewcommand{\t }{\tau }
\renewcommand{\o }{\omega }
\renewcommand{\O }{\Omega }

\newcommand{\norma}[2]{\|#1\|_{\lower 4pt \hbox{$\scriptstyle #2$}}}

\newcommand{\cqd}{{\unskip\nobreak\hfil\penalty50
        \hskip2em\hbox{}\nobreak\hfil\mbox{\rule{1ex}{1ex} \qquad}
        \parfillskip=0pt \finalhyphendemerits=0\par\medskip}}
\newenvironment{pf}{\noindent{\sc Proof}.\enspace}{\rule{2mm}{2mm}\medskip}

\newtheorem{Theorem}{Theorem}[section]
\newtheorem{Corollary}[Theorem]{Corollary}
\newtheorem{Definition}[Theorem]{Definition}
\newtheorem{Lemma}[Theorem]{Lemma}
\newtheorem{Proposition}[Theorem]{Proposition}

\newtheorem{remark}[Theorem]{Remark}
\begin{document}
\title[Porous media equation ]{Some existence and regularity results for porous media and fast diffusion equations with a gradient term}
\author[B. Abdellaoui,  I. Peral, M. Walias]{Boumediene Abdellaoui,  Ireneo Peral and Magdalena Walias}\thanks{Work partially supported by
project MTM2010-18128, MICINN, Spain}

\address{Boumediene Abdellaoui
\hfill \break\indent D\'epartement de Math\'ematiques, Universit\'e Aboubekr Belka\"{\i}d, Tlemcen, \hfill\break\indent Tlemcen 13000,
Algeria.} \email{{\tt boumediene.abdellaoui@uam.es} }

\address{Ireneo Peral
\hfill \break\indent Departamento de Matem{\'a}ticas, U. Aut\'{o}noma de Madrid, \hfill\break\indent 28049 Madrid, Spain.} \email{{\tt
ireneo.peral@uam.es}}

\address{Magdalena Walias
\hfill \break\indent Departamento de Matem{\'a}ticas, U. Aut\'{o}noma de Madrid, \hfill\break\indent 28049 Madrid, Spain.} \email{{\tt
magdalena.walias@uam.es} }

\footnotetext{\it MSC2000, American Mathematical Society, 35D05, 35D10, 35J20, 35J25, 35J70.
} \pagestyle{myheadings} \markboth{}{}
\oddsidemargin=0cm \evensidemargin=0cm \textwidth=160mm \textheight=185mm
\parindent=8mm
\frenchspacing
\date{}
\begin{abstract}
In this paper we consider the problem
$$(P)\qquad
\left\{\begin{array}{rclll} u_t-\D u^m&=&|\n u|^q +\,f(x,t),&\quad u\ge 0 \hbox{ in  } \Omega_T\equiv \Omega\times (0,T),\\ u(x,t)&=&0 &\quad
\hbox{  on  } \partial\Omega\times (0,T)\\ u(x,0)&=&u_0(x),&\quad x\in \Omega
\end{array}
\right.
$$
where $\O\subset \ren$, $N\ge 2$, is a bounded regular domain, $1<q\le 2$, and $f\ge 0$, $u_0\ge 0$ are in a suitable class of functions.

We  obtain some results for elliptic-parabolic problems with measure data related to problem $(P)$ that we use to study  the existence of solutions to problem $(P)$ according with the values of the parameters $q$ and $m$.
\end{abstract}

\maketitle

\rightline{\it To the memory of Juan Antonio Aguilar, our dearest friend.}

\section{Introduction}
In this paper we will study the problem,
\begin{equation}\label{problemaorigen}
\left\{\begin{array}{rclll} u_t-\D u^m&=&|\n u|^q +\,f(x,t),&\quad u\ge 0 \hbox{ in  } \Omega_T\equiv \Omega\times (0,T),\\ u(x,t)&=&0 &\quad
\hbox{  on  } \partial\Omega\times (0,T)\\ u(x,0)&=&u_0(x), &\quad x\in \Omega
\end{array}
\right.
\end{equation}
where $\O\subset \ren$, $N\ge 1 $, is a bounded regular domain, $m>0$, $1<q\le 2$, and $f\ge 0$, $u_0\ge 0$ are in a suitable class of
functions. If  $m>1$, problem \eqref{problemaorigen} is a model of growth in a porous medium, see for instance \cite{baren}.

We refer to the fundamental monograph by J.L. Vázquez, \cite{V}, and the references therein for the basic results about {\it Porous Media Equations} (PME)  and {\it
Fast Diffusion Equation} (FDE)  without gradient term. In \cite{BCP} can be seen a optimal existence result for the Cauchy problem for
the homogeneous (PME).

In this work we are  interested in the existence and regularity of solutions to Problem \eqref{problemaorigen} related to
the parameters $q$ and $m$. Some results were obtained in \cite{DGLS}, \cite{DD}, for $q=2$ and bounded data.
One of the main new features of this paper is to study general class of data according to the values of the parameters in the problems. This
study is in some way motivated by the reference \cite{ADPW}, where the stationary problem is analyzed.

A pertinent formal remark is the following:

{\it The equation $u_t-\Delta u^m=\mu$, after the change $v=u^m$  becomes
$$b(v)_t-\Delta v=\mu \hbox{ with  }\,\, b(s)=s^{\frac 1m}.$$
} The last formulation usually is known in the literature as  {\it elliptic-parabolic equation}. References for problems related to these
equations are \cite{AL}, \cite{Benilan1}, \cite{Benilan2} \cite{BlPo1}, \cite{Carri}, \cite{CP} and \cite{Otto} among others.

We will study the elliptic-parabolic problems with $\mu$ a bounded Radon measure, which is the natural data in the application to the
analysis  of  problem \eqref{problemaorigen}.

The  strategy that we follow in this work can be summarized in the following points.
\begin{enumerate}
\item We consider approximated problems that kill the degeneration or the singularity in the principal  part ((PME) or (FDE) respectively)
and truncating the first order term in the right hand side. With respect to these approximated problems, the existence of a solution follows using the well known results obtained in \cite{BMP1}.
Here the natural setting is to find a {\it weak solution} in the sense of Definition \ref{d1} (where is formulated for the corresponding
elliptic-parabolic equation). \item We obtain uniform estimates of the solution of the approximated problems in such a way that the first
order part in the second member is uniformly bounded in $L^1(\Omega_T)$. \item The previous step motivates the study of a problem with
measure data. To have more flexibility in the calculation we consider the corresponding formulation as elliptic-parabolic equation  and
looking for a {\it reachable solution} in the sense of Definition \ref{d3}. One of the new features in this work, is the proof of the
almost everywhere convergence of the gradient of the solutions of the approximated problems. \item The final step is to use the uniform
estimates and the a.e. convergence of the gradients to prove that, up to a subsequence, the second members of the approximated problems
converge strongly  in $L^1_{loc}(\O_T)$. That is, we find a {\it distributional solution}
\end{enumerate}
The organization of the paper is as follows. In Section 2,  we prove the existence of reachable solutions for a class of {\it
elliptic-parabolic} problem with measure data, including the corresponding cases to the (PME) and (FDE). To obtain the existence of
reachable solutions, we show some a priori estimates for the solutions of the truncated problems and the pointwise convergence of the
gradients, that allow us to conclude. In the last part of Section 2, as an application, we find existence of solution for the porous media
and fast diffusion equations with a Radon measure data. A key point is the proof of a.e. convergence of the gradients, which will be used in
the next section. It is worthy to point out that these results improve in some way the ones obtained in \cite{LKK} for $m>1$, and give a
proof for the (FDE) with measure data.

Section 3 deals with the (PME) with a gradient term. We use the strategy describe above. In the first part of the section we get the
existence results for the interval $1<m\le2$. The main result in this case is Theorem \ref{exis}. In the last part of this section we
consider the complementary interval of $m$, that is, $m>2$. In this case we are able to prove the existence of a solution with
$L^1(\Omega_T)$ data. This result is the content of Theorem \ref{3}.

In the last part of the Section 3, we point out the particular behavior of the case $p=2$, $m=2$. Indeed, if the source term
$f=0$, via a change of
 variable, we show that it is equivalent to a (PME) equation with $m=\frac 53$. As a consequence, we obtain in this case the {\it finite
 speed of propagation}
 property, and also a selfsimilar
solution with compact support in each fixed positive time.

Finally, in Section 4 we analyze the fast diffusion equation, i.e. $0<m<1$, with a convenient hypothesis of integrability of the source term.
The main result of this section is Theorem \ref{fast}. We also prove Theorem \ref{extinct1} that gives us the {\it finite extinction
 time} property of regular solution (see Definition \ref{ddd1}) if $0<m<1$ and $q=2$.

We have tried  to write the paper in an almost  self-contained form,   moreover we give precise references for all the points that are not
detailed in the work.

In a forthcoming paper we will obtain some  results of non-uniqueness for quadratic growth therm.

\section{Some results for an {\it elliptic-parabolic} problem with measure data}
We will consider the general problem
\begin{equation}\label{prob1}
     \left\{\begin{array}{rcll}
 \big(b(v)\big)_t-\Delta\,v&=&\mu
 \quad&\hbox{ in } \ \Omega\times(0,T)\,,\\[5pt]
 v(x,t)&=&0 \quad&\hbox{ on  }
  \p \O\times (0,T)\,,\\[5pt]
  b(v(x,0))&=&b(v_0(x))\quad&\hbox{ in  } \O \,,
 \end{array}\right.
\end{equation}
where $b:\re\to \re$ is continuous strictly increasing function such that $b(0)=0$, $b(v_0)\in L^1(\Omega)$, $f\in L^\infty(\O_T)$ and
 $\mu$ is a Radon measure whose total variation is finite in
$\O_T$. We will assume the following hypotheses on $b$:
$$
(B)\left\{ \begin{array}{lll}(B1) \mbox{ There exists $a_1>0$ such that } b(s)\ge C s^{a_1}\mbox{ for }s>>1,\\(B2)  \mbox{ There exists
$a_2<1$ such that }|b'(s)|\le \dfrac{1}{s^{a_2}}\mbox{ for }s<<1,\\ (B3) \mbox{ There exits }a_3\in (\dfrac{N-1}{N}, 1)\mbox{  such that  }\\
\mbox{ Eithr } b'\in \mathcal{C}([0,\infty))\mbox{  and  }|b'(s)|b^{2a_3-1}(s)\le s^{\frac{N+2a_1}{N}-\e}\mbox{  as  }s\to \infty \mbox{  or
} |b'(s)|\le b^{2-2a_3-\e}(s)\mbox{  as  }s\to \infty.
\end{array}
\right.
$$
\begin{remark}
The following  examples of $b$ will be considered in this work.
\begin{enumerate}
\item  $b(s)=(s+\frac 1n)^{\s}-(\frac 1n)^\s$ if  $s\ge 0$,  for some $\s>\dfrac{(N-2)_+}{N}$  and $1\le n< \infty$ ;
    \item $b(s)=\dyle \frac{1}{m}\int_0^s(H^{-1}(\s))^{\frac{1}{m}-1}d\s=\frac{1}{m}\int_0^{H^{-1}(s)}\s^{\frac{1}{m}-1}H'(\s)d\s$,
    where
$$
H(s)= \dfrac 45 s^{\frac 54} \mbox{  if   }m=2,  \quad \dyle H(s)\int_0^s e^{\frac{t^{\frac{2-m}{m}}}{m(2-m)}}dt  \mbox{ if }0<m<2.
$$
\end{enumerate}
\end{remark}
For $\mu\in L^{\infty}(\Omega_T)$  and $b(v_0)\in L^{\infty}(\Omega)$, we  mean a solution in  the sense of the following definition.

\begin{Definition}\label{d1}
  Assume that $\mu\in L^{\infty}((0,T):L^\infty(\Omega))$
   and $b(v_0)\in L^{\infty}(\Omega)$.
  We say that $v$   is a   weak   solution to \eqref{prob1} if
\begin{enumerate}
\item $v\in L^2((0,T));W^{1,2}_0(\Omega))\cap L^\infty((0,T);L^\infty(\Omega))$, \item The function $b(v)\in\mathcal
    C((0,T);L^{q}(\Omega))$ for all $q<\infty$, \item $(b(v))_t\in     L^{2}((0,T);W^{-1,2}(\Omega))$.
 \end{enumerate}
And for every $\phi\in     L^2((0,T);W^{1,2}_0(\Omega))$  the  following identity holds,
\begin{equation}\label{eq1}     \int_0^T\langle b(v)_t,\phi\rangle+ \iint_{\O_T}\nabla v\cdot\nabla\phi=\iint_{\O_T}\mu\,\phi\,.
 \end{equation}
\end{Definition}

The following result is well known .
 \begin{Theorem}\label{t1}
  Assuming $\mu\in L^{\infty}((0,T);L^\infty(\Omega))$
   and $b(v_0)\in L^{\infty}(\Omega)$,
 there exists a unique weak solution to problem \eqref{prob1} in the sense of Definition \ref{d1}.
 \end{Theorem}

The proof of Theorem \ref{t1} follows  closely the argument developed in \cite{AL} and \cite{CP}.

\subsection{Reachable solutions}\label{sub:reach}

We now introduce the notion of  solution for problem \eqref{prob1} natural to consider measure data. For elliptic equations this notion was
introduced in \cite{DM}. We refer to \cite{Da} for the parabolic equation. See also \cite{ADPS} for some particular cases.

In the case of $\mu \in L^1(\Omega_T)$, the renormalized solution is studied in \cite{BlPo1}. We give a complete analysis under the $(B)$
conditions in order to obtain some regularity on spatial gradient of the solutions.

\begin{Definition}\label{d3}\,
  Assume that $\mu$ is a Radon measure whose total variation is finite
  in  $\O_T$ and $b(v_0)\in L^1(\Omega)$.

  We say that $v$ is a reachable solution to \eqref{prob1} if
  \begin{enumerate}
    \item $T_k(v)\in L^2((0,T);W^{1,2}_0(\Omega))$ for all $k>0$.
       \item For all $t>0$ there exist both one-side limits
 $\lim_{\tau\to t^{\pm}}b(v(\cdot,\tau))$ weakly-* in the sense of measures.
    \item $b(v(\cdot,t))\to b(v_0(\cdot))$ weakly-* in the sense of measures as $t\to0$.
    \item There exist three sequences $\{v_n\}_n$ in
    $L^2((0,T);W^{1,2}_0(\Omega))$, $\{h_n\}_n$ in $L^\infty((0,T);L^\infty(\Omega))$ and
    $\{g_n\}_n$ in $L^\infty(\Omega)$ such that if $v_n$ is the weak solution to problem
    \begin{equation}\label{probla}
     \left\{\begin{array}{rcll}
   \big(b(v_n)\big)_t-\Delta\,v_n&=&h_n
   \quad&\hbox{ in } \ \Omega\times(0,T)\,,\\[5pt]
   v_n(x,t)&=&0 &\hbox{ on  }
  \p \O\times (0,T)\,,\\[5pt]
  v_n(x,0)&=& b^{-1}\big(g_n(x)\big)&\hbox{ in  } \O \,,
   \end{array}\right.
   \end{equation}
   then
    \begin{enumerate}
        \item $g_n\to b(v_0)$ in $L^1(\Omega)$.
        \item $h_n{\buildrel *\over \rightharpoonup}\mu$ as measures.
        \item $\nabla v_n\to \nabla v$ strongly in $L^\sigma((0,T);L^\sigma(\Omega))$
        for $\displaystyle 1\le \sigma<\frac{N+2a_1}{N+a_1}$.
        \item The sequence $\{b(v_n)\}_n$ is bounded in
        $L^\infty((0,T);L^1(\Omega))$ and
        $ b(v_n)\to b(v)$ strongly in $L^1(\O_T)$.
    \end{enumerate}
  \end{enumerate}
 \end{Definition}

To prove the existence of a reachable solution, we need the following Lemma whose proof can be obtained by approximation.

\begin{Lemma}\label{lempart}
 Let $v\in L^2(0,T;W_0^{1,2}(\Omega))$ satisfy $b(v)_t\in
 L^{2'}(0,T;W^{-1,2'}(\Omega))$. Assume that
 $\phi(s)\>:\>\R\to\R$ is a locally Lipschitz--continuous function
 such that $\phi(0)=0$. Then, if we define
 \begin{equation}\label{Fi}
    \Phi(s)=\int_0^s b'(\sigma)\,\phi(\sigma)\,d\sigma\,,
 \end{equation}
 the following integration by parts formula holds:
 \begin{equation}\label{intpart}
    \int_{t_1}^{t_2}\langle b(u)_t,\phi(u)\rangle_{W^{-1,2'}(\Omega),
    W_0^{1,2}(\Omega)}\,dt=
    \int_\Omega\Phi(u(x, t_2))\,dx-\int_\Omega\Phi(u(x,
    t_1))\,dx\,,
\end{equation}
 for every $0\le t_1<t_2$.
\end{Lemma}

We will now prove the existence of a reachable solution to Problem \eqref{probla}.

\subsection{Some a priori estimates}\label{subsec:}
Let us consider the following approximating problems:
\begin{equation}\label{eq:pbap}
\left\{
\begin{array}{rcll}
b(v_n)_t-\Delta v_n&=&h_n,\, &(x,t)\in \O_T,\\ v_n(x,t)&=&0,\, &(x,t)\in\partial\Omega\times(0,T),\\ v_n(x,0)&=&b^{-1}(g_n(x)),\,
&x\in\Omega,
\end{array}
\right.
\end{equation}
where $g_n\to b(v_0)$ strongly in $L^1(\Omega)$ and $h_n\to \mu$ in the weak-$*$ sense in $\Omega_T$. The existence of weak solutions to
these problems follows from Theorem \ref{t1}.

Let us begin by proving the next proposition.
\begin{Proposition}\label{stima}
Let $\{v_n\}_n$ be a sequence of solutions of the approximate problems \eqref{eq:pbap}. Then
\begin{enumerate}
\item For each $0<\beta<\frac 12$, the sequence $\{(|v_n|+1)^\beta-1\}_n$ is bounded in $L^2(0,T; W^{1,2}_0(\Om))$ \item The sequence
    $\{b(v_n)\}_n$ is bounded in the space $L^\infty(0,T; L^1(\Om))$ and $\{(b(v_n))_t\}_n$ is bounded in
    $L^1(\O_T)+L^\sigma(0,T;W^{-1,\sigma}(\O)), \hbox{ for some } \sigma>1$.
\end{enumerate}
Moreover,
\begin{equation}\label{stima3m}
 \iint_{\O_T}\frac{|\n v_n|^2}{(1+|v_n|)^{\a+1}} \le C\qquad\hbox{for all }
\a>0.
\end{equation}
Furthermore, the sequence $\{|\nabla v_n|\}_n$ is bounded in the Marcinkiewicz spaces $M^{q}(\O_T)$ where $q=\frac {(N+2a_1)}{N+a_1}$ and
$\{v_n\}_n$ is bounded in the space $M^{\sigma}(\O_T)$ where $\sigma=\frac {N+2a_1}N$.
\end{Proposition}

\begin{pf} Take $\phi(v_n)\equiv T_k(v_n)$ in Lemma \ref{lempart}, with $k>0$, as test function in
the weak formulation of \eqref{eq:pbap}, then
 \begin{equation}\label{ine1}
    \int_\O\Phi_k(v_n(x, t))\,dx-\int_\O\Phi_k(b^{-1}(g_n(x)))\,dx+
   \iint_{\O_T}|\nabla T_k(v_n)|^2\le k|\mu|(\O_T)\,,
\end{equation}
where $\Phi_k(s)=\int_0^sT_k(\sigma)\, b^\prime(\sigma)\,d\sigma$ which is a nonnegative function. Since $\Phi_k(s)\le k|b(s)|$, it follows
from \eqref{ine1}  that
 \begin{equation}\label{ine2}
 \int_\O|\Phi_k(v_n(x, t))|\,dx+
    \iint_{\O_T}|\nabla T_k(v_n)|^2\le
    ck\iint_{\O_T}(1+|b(v_n)|)+ck\,.
 \end{equation}
 Now, dropping a nonnegative term, dividing by $k$ and letting $k$
 go to $0$, it yields
  $$
  \int_\O|b(v_n(x, t))|\,dx\le
    c \iint_{\O_T}|b(v_n)|+c\,.
    $$
 Thus, Gronwall's Lemma implies that
 \begin{equation}
    \sup_{t\in[0,T]}\int_\O|b(v_n(x, t))|\,dx\le C.\label{stima1}
    \end{equation}
    Moreover, going back to \eqref{ine2} we get
    \begin{equation}
     \iint_{\O_T}|\nabla T_k(v_n)|^2\le  Ck\,.\label{stima2}
 \end{equation}

In order to prove the next estimate, we define $\theta(s)=\big(1-(1+|s|)^{-\a}\big)$, with $0<\alpha$, and we take
 $\theta(v_n)$ as test function in the approximating problems. There
 results,
\begin{equation}\label{eq:0}
\io \Theta(v_n(T)) dx-\io \Theta(b^{-1}(g_n)) dx+\alpha \iint_{\O_T}\frac{|\n v_n|^2}{(1+|v_n|)^{\a+1}} \le ||\mu||,
\end{equation}
 where $\Theta(s)=\int_0^s b'(\sigma)\,\theta(\sigma)\,d\sigma$. Hence
  \begin{equation*}
   \iint_{\O_T}\frac{|\n v_n|^2}{(1+|v_n|)^{\a+1}}\le C\mbox{  for all  }\a>0
 \end{equation*}
 and then $\{(|v_n|+1)^\beta-1\}_n$ is bounded in $L^2(0,T;
W^{1,2}_0(\Om))$ for all $0<\beta<\frac 12$.

To get the estimates on the Marcinkiewicz spaces, we follow closely the arguments in \cite{B-V}, see also \cite{AMST}.

From  hypothesis  $(B1)$, we obtain that $b(s)\ge C_1 s^{a_1} -C\mbox{ for all }s\ge 0$. By  using again \eqref{stima1} we get the existence
of a   positive constant $C$ such that,
  \begin{equation}\label{stMarc1}
|\{x\in\Omega\>:\>|v_n(x, t)|>k\}|\le \frac {C}{k^{a_1}}\quad\hbox{for almost all
  }t\in[0,T]\,,\hbox{ all }k>0\hbox{ and all }n\in\N.
  \end{equation}
 Thus, by Sobolev's inequality and \eqref{stima2},
  \begin{multline}\label{stMarc2}
  \int_0^T\big(|\{x\in\Omega\>:\>|v_n(x, t)|\ge k\}|\big)^{2/2^*}\,dt\le
  \int_0^T\bigg(\frac{\|T_k(v_n(x, t))\|_{2^*}^{2^*}}{k^{2^*}}\bigg)^{2/2^*}\,dt\\
    \le C\int_0^T\frac{\|\nabla
  T_k(v_n(x, t))\|_{2}^{2}}{k^{2}}\,dt\le \frac C{k}\quad
  \hbox{ for all }k>0\hbox{ and all }n\in\N\,.
  \end{multline}
  Therefore, combining \eqref{stMarc2} and \eqref{stMarc1} we
  obtain, for all $k>0$ and all $n\in\N$,
  \begin{equation}\label{stMarc3}
\begin{array}{ll}
  &|\{(x, t)\in \O_T\>:\>|v_n(x, t)|\ge k\}|\\
   &\dyle=\int_0^T\big(|\{x\in\Omega\>:\>|v_n(x, t)|\ge k\}|\big)^{1-(2/2^*)}
  \big(|\{x\in\Omega\>:\>|v_n(x, t)|\ge k\}|\big)^{2/2^*}\,dt\le \frac C{k^{\frac{N+2a_1}{N}}}.
  \end{array}
  \end{equation}
  A similar estimate for the gradients is now easy to obtain.
  Indeed, for every $h,k>0$, we have
$$
\begin{array}{ll}
&|\{(x, t)\in \O_T\>:\>|\nabla v_n(x, t)|\ge h\}|\\ &\le |\{(x, t)\in \O_T\>:\>|v_n(x, t)|\ge k\}|+|\{(x, t)\in \O_T\>:\>|\nabla T_k(v_n(x,
t))|\ge h\}| \\ &\le \dyle  \frac C{k^{\frac{N+2a_1}{N}}}+ \iint_{\O_T}\frac{|\nabla
  T_k(v_n)|^2}{h^2}\le
  \frac C{k^{\frac{N+2a_1}{N}}}+\frac {Ck}{h^2}\,.
  \end{array}
  $$
 Minimizing in  $k$, we obtain that for  $k=h^{\frac{N}{N+a_1}}$,
 \begin{equation}\label{stMarc4}
 |\{(x, t)\in \O_T\>:\>|\nabla v_n(x, t)|\ge h\}|\le \frac
 C{h^{\frac{N+2a_1}{N+a_1}}}\,
  \hbox{ for all }h>0\hbox{ and all }n\in\N\,.
 \end{equation}
Hence,
  \begin{gather}
     \iint_{\O_T}|v_n|^\rho\le C\quad\hbox{for all } 0<\rho<\frac{N+2a_1}{N}\label{stima4}\\
     \iint_{\O_T}|\nabla v_n|^r\le C\quad\hbox{for all } 0<r<\frac{N+2a_1}{N+a_1}.\label{stima5}
 \end{gather}
From  equations \eqref{eq:pbap} we obtain,
$\{\big(b(v_n)\big)_t\}_n \hbox{ is bounded in }   L^1(\O_T)+L^\sigma(0,T;W^{-1,\sigma}(\O))$ for some $ \sigma\ge \frac{N+2a_1}{a_1}.$
Consider $\varrho$  such that
 \begin{equation}\label{euro000}
-\D \varrho =1\mbox{  in   }\O,   \,\varrho\in \sob(\O).
\end{equation}
We claim that for all $0<\d<\min\{1, \frac{1}{a_2}\}$,
 \begin{equation}\label{euro}
 \int_0^{T}\int_{\O}\dfrac{|\nabla v_n|^2}{(v_n+\frac 1n)^{1+\d}}\varrho\le
 C\mbox{  uniformly  in  }n.
\end{equation}
To prove the claim we define $K(s)=\dyle\int_0^s \frac{b'(\s)}{(\s+\frac 1n)^\d} d\s.$ Using $(B2)$ we get easily that $K(v_n)$ is well
defined, $K(0)=0$ and
 \begin{equation}
    \sup_{t\in[0,T]}\int_\O K(v_n(x, t))\varrho\,dx\le C.\label{euro2},
    \end{equation}
Using $\dfrac{\varrho}{(v_n+\frac 1n)^{\d}}$ as a test function in \eqref{eq:pbap}, it follows that
\begin{eqnarray*}\label{euro00}
&\dyle \io K(v_n(x,T))\varrho dx +\frac{1}{1-\d} \iint_{\O_T}\Big[(v_n+\frac 1n)^{1-\d}-(\frac 1n)^{1-\d}\Big](-\D \varrho )=\\ &\dyle
\iint_{\O_T}\dfrac{|\nabla v_n|^2}{(v_n+\frac 1n)^{1+\d}}\varrho + \iint_{\O_T}\dfrac{h_n}{(v_n+\frac 1n)^{1+\d}}\varrho+\io
K(v_n(x,0))\varrho dx.
\end{eqnarray*}
Since $0<\d<1$, then using \eqref{stima4}, we get
$$
\frac{1}{1-\d} \iint_{\O_T}\Big[(v_n+\frac 1n)^{1-\d}-(\frac 1n)^{1-\d}\Big]\le C\mbox{  uniformely  in   }n.
$$
 As a consequence, and by \eqref{euro2},
    \begin{equation}
 \iint_{\O_T}\dfrac{|\nabla v_n|^2}{(v_n+\frac 1n)^{1+\d}}\varrho\le
C_1\mbox{  uniformly  in  }n. \label{euro3}
 \end{equation}
In particular the claim follows.

In the same way, using $1-\dfrac{1}{(b(v_n)+1)^{\d}}$, where $0<\d$, as a test function in \eqref{eq:pbap}, it follows that
\begin{equation}
\sup_{t\in[0,T]}\int_\O b(v_n(x, t))\,dx+ \int_0^{T_1}\io \dfrac{b'(v_n)|\nabla v_n|^2}{(b(v_n)+1)^{1+\d}}\le C\mbox{ uniformly in  }n.
\label{euro5}
 \end{equation}

 To finish we have just to prove that
 \begin{equation}\label{strong}
 b(v_n)\to b(v)\mbox{  strongly  in } L^1(\O_T).
 \end{equation}
 Notice that $||b(v_n)||_{L^1(\O_T)}\le C$ and $b(v_n)\to b(v)$
 e.a in $\O_T$.

If $b(s)\le C s^{\frac{N+2a_1}{N}-\e}$, for some $\e>0$, as $s\to \infty$, then using \eqref{stima4} and by Vitali's lemma we reach the
strong convergence in \eqref{strong}.

Assume that  the condition $B_3$ holds, let $w_n=(b(v_n)+1)^\beta$, where $\beta\le 1$ to be chosen later, then
$||w_n||_{L^\infty(0,T;L^1(\O))}\le C$ and $\n w_n =\beta b'(v_n)(b(v_n)+1)^{\beta-1}\n v_n$. Thus, for $\d>0$, we have
\begin{eqnarray*}
\dyle \iint_{\O_T}|\n w_n| &= &\beta \iint_{\O_T}b'(v_n)(b(v_n)+1)^{\beta-1}|\n v_n|=\beta
\iint_{\O_T}b'(v_n)(b(v_n)+1)^{\frac{2\beta+\d-1}{2}}\dfrac{|\n v_n|}{(b(v_n)+1)^{\frac{\d+1}{2}}}\\ &\le & \dyle \dfrac 12 \beta
\iint_{\O_T}\dfrac{b'(v_n)|\n v_n|^2}{(b(v_n)+1)^{\d+1}} + \dfrac 12 \beta \iint_{\O_T}b'(v_n)(b(v_n)+1)^{2\b+\d-1}.
\end{eqnarray*}
From \eqref{euro5}, we obtain that
$$
\iint_{\O_T}\dfrac{b'(v_n)|\n v_n|^2}{(b(v_n)+1)^{\d+1}}\le C.
$$
Now, if the first condition in $B_3$ holds, then choosing $\d$ small enough and $\beta=a_3$, we reach that
$$
b'(s)(b(s)+1)^{2\b+\d-1}\le C_1 s^{\frac{N+2a_1}{N}-\e}\mbox{  as  }s\to \infty \mbox{  for some  }\e>0.
$$
Thus
$$
\dfrac{\beta}2  \iint_{\O_T}b'(v_n)(b(v_n)+1)^{2\b+\d-1}\le C\mbox{  uniformly in  }n.
$$
If the second condition in in $B_3$ holds, with the same chose of $\beta$ and for for $\d$ small,
$$
b'(s)(b(s)+1)^{2\b+\d-1}\le b^{1-\e}(s)\mbox{  as   }s\to \infty,
$$
and hence
$$
\dfrac {\beta}2\iint_{\O_T}b'(v_n)(b(v_n)+1)^{2\b+\d-1}\le C\mbox{  uniformly in  }n.
$$

Therefore, we conclude that $||w_n||_{L^1(0,T;W^{1,1}0(\O))}\le C$ for all $n$. Hence, using the Gagliardo-Nirenberg inequality we conclude
that $||w_n||_{L^{\frac{N}{N-1}}(\O_T)}\le C$.

Since $\beta \frac{N}{N-1}>1$, then $\{b(v_n)\}_n$ is bounded in $L^{1+\e}(\O_T)$ for some $\e>0$.
 Hence, using Vitali's lemma, we obtain that the sequence
 $\{b(v_n)\}_n$ is compact in $L^1(\O_T)$ and so we may
 extract a subsequence (also labelled by $n$) such that
 $b(v_n)\to b(v)$ strongly in $L^1(\O_T)$.
\end{pf}
\begin{remark}
If $b(s)=(s+\frac 1n)^{\s}-(\frac 1n)^\s$, then
\begin{enumerate}
\item If  $\s\le \frac{N}{N-1}$, we can prove that the sequence $\{\big(b(v_n)\big)\}_n$ is bounded in $L^1(0,T;W^{1,1}_{loc}(\O))$. This
    follows using estimates \eqref{euro5}. \item If  $\s\le \frac{N}{N-2}$, we can prove that the sequence $\{\big(b(v_n)\big)\}_n$
    converge strongly in $L^1(\O_T)$. This follows using \eqref{stima4}.
\end{enumerate}
\end{remark}

\subsection{Pointwise convergence of the gradients}
In this subsection we prove that up to a subsequence
$$\nabla v_n\to \n v,\,  \hbox{ a.e. in }    \O_T \hbox{ as  }  n\to \infty,.$$
Hence we will obtain that the sequence $\{v_n\}_n$ satisfies condition (4) (c) in Definition \ref{d3}.

\begin{Proposition}\label{convgr}
Consider $\{v_n\}_n$, the solution of the approximated problems \eqref{eq:pbap}. Then, up to subsequence,
\begin{equation}\label{conv2}
    \nabla T_{k}(v_n)\to\nabla T_{k}(v)\qquad\hbox{almost everywhere
    in } \O_T\,.
\end{equation}
As a consequence,  $\nabla v_n\to\nabla v$ almost everywhere in $\O_T.$
\end{Proposition}

\begin{pf} We recall the
time--regularization of functions due to Landes and Mustonen (see \cite{L}, \cite{LM}). Consider $w$ such that $T_k(w)\in
L^2(0,T:\sob(\Omega))\cap \mathcal{C}([0,T]:L^2(\Omega))$. For every $\nu\in \N$, we define $(T_k w)_\nu$ as the solution of the Cauchy
problem
\begin{equation}\label{regg}
\left\{
\begin{array}{ll}
\displaystyle {1\over\nu}[(T_kw)_\nu]_t+(T_kw)_\nu = T_{k}w;\\
\\
(T_kw)_\nu(0)=T_k(w_0).\\
\end{array}
\right.
\end{equation}
Then, one has,
$
   (T_k w)_\nu \in L^2(0,T;W^{1,2}_0(\Om))$, $((T_kw)_\nu)_t \in L^2(0,T;W^{-1,2}_0(\Om))$,
$$ \Vert
(T_kw)_\nu\Vert_{L^\infty(\O_T)}\le\Vert T_kw\Vert_{L^\infty(\O_T)} \leq k, $$ and as $\nu$ goes to infinity, $(T_kw)_\nu \to T_k w$ strongly
in $L^2(0,T;W^{1,2}_0(\Om)).$

From now on, we will denote by $\omega(n, \nu, \varepsilon)$ any quantity such that
  $$
  \lim_{\varepsilon\to0^+}\limsup_{\nu\to\infty}\limsup_{n\to\infty}
  |\omega(n, \nu, \varepsilon)|=0\,.
  $$
Taking $T_\varepsilon\big(b(v_n)-(T_k(b(v))_\nu\big)$ as test function in \eqref{eq:pbap}, we obtain
\begin{equation*}
\begin{array}{lll}
& &\dyle\int_{0}^{T} \langle(b(v_n))_t\,,\,T_\varepsilon\big(b(v_n)-(T_k(b(v)))_\nu\big)\big)\rangle\,dt+
 \iint_{\O_T}\nabla v_n\cdot \nabla
 T_\varepsilon\big(b(v_n)-(T_k(b(v))_\nu\big)\big)\\
 &=&\dyle\iint_{\O_T}h_n \, T_\varepsilon\big(b(v_n)-(T_k(b(v)))_\nu\big)
\end{array}
\end{equation*}
  We will analyze the integrals which appear in the previous
 equality.
For simplicity of typing we set $w_n=b(v_n)$,
 $w=b(v)$ and $\psi =b^{-1}$.

It is clear that
$$
\iint_{\O_T}h_n \, T_\varepsilon\big(b(v_n)-(T_k(b(v)))_\nu\big)\le C\e.
$$
Notice that
$$\int_{0}^{T}
\langle(b(v_n))_t\,,\,T_\varepsilon\big(b(v_n)-(T_k(b(v)))_\nu\big)\big)\rangle\,dt=\int_{0}^{T}
\langle((w_n)_t\,,\,T_\varepsilon\big(w_n-(T_k(w))_\nu\big)\big)\rangle\,dt.
$$
Then using the same arguments as in \cite{BDGO} (see also (3.37) in \cite{ADPS}) we reach that
\begin{equation}\label{stimA}
\int_{0}^{T} \langle(b(v_n))_t\,,\,T_\varepsilon\big(b(v_n)-(T_k(b(v)))_\nu\big)\big)\rangle\,dt\ge\omega(n,\nu, \eps)\,.
\end{equation}
We set
$$
I= \iint_{\O_T}\nabla v_n\cdot \nabla T_\varepsilon\big(b(v_n)-(T_k(b(v))_\nu\big)\big).
$$
{\it Claim.}  It holds that
  \begin{equation}\label{stimB}
 I\geq \iint_{\{|w_n|\le k\}}\psi'(T_k(w_n))\nabla T_k(w_n)
  \cdot \nabla T_\varepsilon(T_k(w_n) - (T_k w)_\nu)
  +\omega^\varepsilon(n,\nu)\,.
  \end{equation}
Indeed,
  \begin{equation*}
  \begin{array}{lll}
  I &=& \dyle \iint_{\{|w_n|\le k\}} \psi'(w_n)\nabla T_k(w_n)
  \cdot \nabla T_\varepsilon\big(T_k(w_n) - (T_k w)_\nu\big)\\ &+&\dyle \iint_{\{|w_n|> k\}} \psi'(w_n)\nabla w_n
  \cdot \nabla T_\varepsilon\big(w_n - (T_k w)_\nu\big)
  \\
  &\geq &\dyle \iint_{\{|w_n|\le k\}} \psi'(T_k(w_n))\nabla T_k(w_n)
  \cdot \nabla T_\varepsilon\big(T_k(w_n) - (T_k w)_\nu\big)
 \\
  &-& \dyle \iint_{\{|w_n|> k\,,\ |w_n - (T_k w)_\nu|\leq \eps \}} \psi'(w_n)\nabla w_n
  \cdot \nabla  (T_k w)_\nu\,.
  \end{array}
  \end{equation*}
Since $\|(T_k w)_\nu\|_\infty\leq k$, the last integrand is different from zero only in the set where $|w_n|\leq k+\eps$, therefore the last
integrand is bounded by
  \begin{equation*}
  \begin{split}
   &c_1(k)\,\left[\iint_{\O_T} |\nabla T_{k+\eps}w_n|^2\right]^{\frac{1}{2}}\,
   \left[\iint_{\O_T} |\nabla (T_{k}w)_\nu|^2 \car_{\{|w_n|>k\}}\right]^{\frac{1}{2}}\\
   &\qquad \le c_2(k,\eps)\,\left[\iint_{\O_T} |\nabla T_{k}w|^2 \car_{\{|w|>k\}}\right]^{\frac{1}{2}} + \omega^{\eps}(n, \nu)=
   \omega^{\eps}(n, \nu)\,.\\
     \end{split}
  \end{equation*}
Here we have used the a.e. convergence of $\car_{\{|w_n|>k\}}$ to $\car_{\{|w|>k\}}$ as $n\to+\infty$, which holds for all $k$ (see for
instance \cite{BDGO}).

Therefore, putting together the above estimates, we reach that
\begin{equation}\label{conv5}
\begin{array}{ll}
&\dyle\iint_{\O_T} \psi'(T_k(w_n))\nabla T_k(w_n)
  \cdot \nabla T_\varepsilon\big(T_k w_n - (T_k w)_\nu\big)\\&=
 \dyle\iint_{\{|w_n|\le k\}} \psi'(w_n)\nabla T_k(w_n)
  \cdot \nabla T_\varepsilon\big(w_n - (T_k w)_\nu\big)\le\omega(n,\nu, \varepsilon)\,.
\end{array}
\end{equation}
Thus
  \begin{equation}\label{conv6}
\iint_{\O_T} \psi'(T_k(w_n))\nabla T_k(w_n)
  \cdot \nabla \big(T_k w_n - (T_k w)_\nu\big)\chi_{{\{|T_k(w_n)-(T_kw)_\nu|\le\varepsilon\}} }\le\omega(n,\nu, \varepsilon)\,.
\end{equation}
 Hence it follows that
\begin{multline*}
    \iint_{\O_T} \psi'(T_k(w_n))\nabla T_k(w_n)
  \cdot \nabla \big(T_k w_n - T_k w\big)\chi_{\{|T_k(w_n)-(T_kw)_\nu|\le\varepsilon\}}\\
    =\iint_{\O_T} \psi'(T_k(w_n))\nabla T_k(w_n)
  \cdot \nabla \big(T_k w_n - (T_k   w)_\nu\big)\chi_{\{|T_k(w_n)-(T_kw)_\nu|\le\varepsilon\}}\\
  +\iint_{\O_T} \psi'(T_k(w_n))\nabla T_k(w_n) \cdot \nabla \big((T_k w)_\nu-T_kw\big)\chi_{\{|T_k(w_n)-(T_kw)_\nu|\le\varepsilon\}}\\
 = \iint_{\O_T} \psi'(T_k(w_n))\nabla T_k(w_n)
  \cdot \nabla \big(T_k w_n - (T_k   w)_\nu\big)\chi_{\{|T_k(w_n)-(T_kw)_\nu|\le\varepsilon\}}   +\omega(n,\nu, \varepsilon).
\end{multline*}

  Therefore, by \eqref{conv6} we  obtain
    \begin{equation}\label{conv7}
\iint_{\O_T} \psi'(T_k(w_n))\nabla T_k(w_n)
  \cdot \nabla \big(T_k w_n - T_k w\big)\chi_{\{|T_k(w_n)-(T_kw)_\nu|\le\varepsilon\}} \le\omega(n,\nu, \varepsilon)\,.
\end{equation}

On the other hand, we have that
 \begin{equation*}
 \iint_{\O_T} \psi'(T_k(w_n))\nabla  T_kw
  \cdot \nabla \big(T_k w_n - T_k
  w\big)\chi_{\{|T_k(w_n)-(T_kw)_\nu|\le\varepsilon\}}=\omega^{\nu,
  \varepsilon}(n)\,.
\end{equation*}
Indeed, $\psi'(T_k(w_n))\nabla  T_kw\to \psi'(T_k(w))\nabla  T_kw$ strongly in $L^{2}(\O_T)$ and
  $\nabla T_k w_n \rightharpoonup \nabla T_k w$ weakly in $L^{2}(\O_T)$. Hence, we deduce from \eqref{conv7} that
      \begin{equation}\label{conv7.5}
  \iint\limits_{\{|T_k(w_n)-(T_kw)_\nu|\le\varepsilon\}} \!\! \psi'(T_k(w_n))|\nabla T_k(w_n)-\nabla  T_kw|^2\le\omega(n,\nu,
  \varepsilon)\,.
\end{equation}
Denoting $\Psi_{n,k}=\psi'(T_k(w_n))|\nabla T_k(w_n))-\nabla T_kw|^2$, then
  \begin{equation}\label{conv8}
    \iint_{\{|T_k(w_n)-(T_kw)_\nu|\le\varepsilon\}}\Psi_{n,k}\le\omega(n,\nu, \varepsilon)\,.
\end{equation}
Since
 \begin{equation}\label{cond}
    \chi_{\{|T_k(w_n)-(T_kw)_\nu|>\varepsilon\}}\to
    \chi_{\{|T_k(w)-(T_kw)_\nu|>\varepsilon\}}\qquad\hbox{strongly in
    } L^\rho(\O_T)\,,\forall \rho\ge1\,.
\end{equation}
then using H\"older inequality and the fact that $\Psi_{n,k}$ is bounded in $L^1(\O_T)$,  we deduce that

$\displaystyle \iint_{\O_T}\Psi_{n,k}^\theta=\omega(n, \nu, \varepsilon)$ for all $\theta<1$. Since $\psi'(T_k(s))\ge c(k)(\psi'(T_k(s)))^2$, it follows
that
 \begin{equation}\label{TLM}
 \lim_{n\to\infty}\iint_{\O_T}\Big[(\psi'(T_k(w_n)))^2|\nabla
T_k(w_n)-\nabla  T_kw|^2\Big]^\theta=0\,,
\end{equation}
 Hence we get $\psi'(T_k(w_n))\n T_k(w_n)\to \psi'(T_k(w))\n
 T_k(w)$ e.a in
 $\O$ and then
 $$\n \psi(T_k(w_n))\to \n \psi(T_k(w_n))\quad\hbox{ a.e. in }\O.$$

  Using the fact that $\psi$ is a strictly monotone function we reach
$$\n T_k(w_n)\to \n T_k(w_n)\quad\hbox{  a.e. in } \O$$
 and then $\n T_k(v_n)\to \n T_k(v_n)$ a.e in  $\O$.
\end{pf}

\begin{Corollary}\label{convgr2}

\

\begin{enumerate}
\item Let $\{v_n\}_n$ be a sequence of solutions of \eqref{eq:pbap}, then $\nabla v_n\to \nabla v$ strongly in $L^\sigma(\Omega_T)$ for
    all $1\le \sigma<\frac{N+2a_1}{N+a_1}$.

\item From estimate \eqref{TLM}, we reach
\begin{equation}\label{TLMM}
\n T_k(v_n)\to \n T_k(v)\mbox{  strongly  in   }L^\s(\O_T)\mbox{ for all  }\s<2.
\end{equation}
\end{enumerate}
\end{Corollary}
\begin{Corollary}\label{distr}
    The
    following equality holds
    \begin{equation}\label{distn}
    -\int_\Omega b(v_0(x))\Phi(x,0)\,dx-\iint_{\O_T}b(v)\Phi_t
    +\iint_{\O_T}\nabla v\cdot\nabla\Phi=\iint_{\O_T}
    \Phi \,d\mu\,,
   \end{equation}
   for every $\Phi\in C^\infty(\overline \O_T)$, with $\Phi(\cdot,t)\in C_0(\Omega)$ for all $t\in(0,T)$ and $\Phi(x,T)=0$ for all
   $x\in\Omega$.
\end{Corollary}
\subsection{Application to the porous media and fast diffusion  equations with a Radon measure.}\label{sec:R}
The results in the above subsection allow us to consider the problem
\begin{equation}\label{prob14}
     \left\{\begin{array}{rcll}
 u_t-\Delta\,u^m&=&\mu
 \quad&\hbox{ in } \ \Omega\times(0,T)\,,\\[5pt]
 u(x,t)&=&0 \quad&\hbox{ on  }
  \p \O\times (0,T)\,,\\[5pt]
  u(x,0)&=&u_0(x)\quad&\hbox{ in  } \O \,,
 \end{array}\right.
\end{equation}
with $m>\dfrac{(N-2)_+}{N}$, $u_0\in L^1(\O)$ and $\mu$ is   a Radon measure whose total variation is finite in $\O_T$.

In the case of the porous media equation, i.e., $m>1$, the existence results are obtained in \cite{LKK} by using some result in \cite{KL}.
Here we extend the results to the interval $\dfrac{(N-2)_+}{N}<m<10$, by proving, moreover, the a.e convergence of the gradients of the
truncated problems to the gradient of the solution of problem \eqref{prob14}. Our approach is  alternative by using the {\it
elliptic-parabolic}  framework.

We will consider the approximated form
\begin{equation}\label{prob141}
     \left\{\begin{array}{rcll}
 u_{nt}-\text{div}(m(u_n+\frac 1n )^{m-1}\n u_n)&=&h_n
 \quad&\hbox{ in } \ \Omega\times(0,T)\,,\\[5pt]
 u_n(x,t)&=&0 \quad&\hbox{ on  }
  \p \O\times (0,T)\,,\\[5pt]
  u_n(x,0)&=&T_n(u_0)\quad&\hbox{ in  } \O.
 \end{array}\right.
\end{equation}
The main goal of this subsection is to show compactness results for the sequences $\{|\n u_n|\}_n$ and $\{T_k(u_n)\}_n$ including the case
$\dfrac{(N-2)_+}{N}<m<1$.

Define $v_n\equiv (u_n+\frac 1n )^{m}-(\frac 1n)^m$, then $v_n$ solves
\begin{equation}\label{eq:pbap1}
\left\{
\begin{array}{rcll}
(b(v_n))_t-\Delta v_n&=&h_n,\, &(x,t)\in \O_T,\\ v_n(x,t)&=&0,\, &(x,t)\in\partial\Omega\times(0,T),\\
v_n(x,0)&=&\varphi^{-1}(T_n(u_0(x))),\, &x\in\Omega,
\end{array}
\right.
\end{equation}
where $b(s)=(s+(\frac 1n)^m)^{\frac 1m}-\frac 1n$.

\begin{Theorem}\label{last00}
Consider $v_n$, the solution to \eqref{eq:pbap1} and $u_n=b(v_n)$ that solves\eqref{prob141}. There exists a measurable function $u$ such
that $u^m\in L^r(0,T;W^{1,r}_0(\Om))$ for all $r<1+\frac{1}{Nm+1}$, and, up to a subsequence,
\begin{enumerate}
\item $\n u_n\to \n u$ e.a in $\O_T$ and then $\n v_n\to \n v$ e.a in $\O_T$ where $v=u^m$. \item $T_k(v_n)\to T_k(v)$ strongly in
    $L^\s(0,T;W^{1,\s}_0(\O))$ for all $k>0$ and for all $\s<2$.
\end{enumerate}
\end{Theorem}
\pf  Since $b$ verifies the condition $(B)$ with $a_1=\frac 1m$, then by the results of Corollary \ref{convgr2}, we obtain  that
\begin{enumerate}
\item $\nabla v_n\to \nabla v$ strongly in the space $L^\sigma((0,T);L^\sigma(\Omega))$ for all $1\le \sigma<1+\frac{1}{Nm+1}$ and \item
    $\n T_k(v_n)\to \n T_k(v)$  strongly  in   $L^\s(\O_T)$ for all  $\s<2$.
\end{enumerate}
Let $\psi\in \mathcal{C}^\infty_0(\O_T)$ be such that $\psi\ge 0$ in $\O_T$, we claim that \begin{equation}\label{UAM1} \Big\{\dfrac{|\n
T_k(u_n)|}{(u_n+\frac 1n)^\theta}\psi\Big\}_n\mbox{  is uniformely bounded in } L^1(\O_T)\mbox{  for some }\theta>0.
\end{equation}
To prove the claim we will consider separately the cases: $0<m<1$ and $m>1$.

Let begin by the case $0<m<1$. Using $T_k(u_n)$ as a test function in \eqref{prob14}, we get
$$
\iint_{\O_T}(u_n+\frac 1n)^{m-1}|\n T_k(u_n)|^2\le C.
$$
Thus the claim follows with $\theta =1-m>0$.

We deal now with the case $m>1$.   Using $\dfrac{\psi}{(v_n+(\dfrac 1n)^m)^\d}$, where $\d<\dfrac 1m$, as a test function in \eqref{eq:pbap1}
we obtain that
\begin{equation}\label{FGH}
\iint_{\O_T} \dfrac{|\n v_n|^2}{(v_n+(\frac 1n)^m)^{\d+1}}\psi\le C.
\end{equation}
Notice that
$$
\dfrac{|\n T_k(u_n)|}{(u_n+\frac 1n)^\theta}\psi=\dfrac{|\n T_{k'} v_n|}{(v_n+(\frac 1n)^m)^{1-\frac{1-\theta}{m}}}\psi .
$$
Now using \eqref{FGH} we conclude that $\Big\{\dfrac{|\n T_k(u_n)|}{(u_n+\frac 1n)^\theta}\psi\Big\}_n$ is bounded in $L^1(\O_T)$ for all
$\theta<1$. Then, the claim is proved.

By similar computation as above, taking $\dfrac{\psi}{(u_n+\frac 1n)^\d}$ as a test function in \eqref{prob141} with $\d<1$, we reach that
\begin{equation}\label{UAMM}
\d\iint_{\O_T} \dfrac{|\n u_n|^2}{(u_n+\frac 1n)^{\d+1-m}}\psi+\iint_{\O_T} \dfrac{h_n}{(u_n+\frac 1n)^\d}\psi\le C\mbox{ for all }n.
\end{equation}
We will prove  the point (1) in the Theorem, that is $\n u_n\to \n u$ a.e in $\O_T$  where $u=v^{\frac 1m}$.

To prove this assertion   it is sufficient to show that, for some $s\in (0, 1)$,
\begin{equation}\label{S}
\dyle\io|\n T_k(u_n)-\n T_k(u)|^s \psi \to 0\mbox{  as   }n\to \infty.
\end{equation}
Consider $A\equiv\{(x,t)\in\O_T:\, u(x,t)=0\}\equiv \{(x,t)\in\O_T:\, v(x,t)=0\}$. Notice that $\nabla u_n\to \nabla u$ in $\O_T\setminus A$.
Since $\n T_k(u)=0$ and $\n T_k(v)=0$ in $A$, then
\begin{equation}\label{RR}
\dyle \iint_{\O_T}|\n T_k(u_n)-\n T_k(u)|^s\psi = \dyle \iint_{\{u=0\}} |\n T_k(u_n)|^s \psi +\dyle \iint_{\{u>0\}} |\n T_k(u_n)-\n T_k(u)|^s
\psi
\end{equation}
Since $s<1$ and $\nabla u^m_n\to \nabla u^m$ strongly in the space $L^\sigma((0,T);L^\sigma(\Omega))$ for all $1\le \sigma<1+\frac{1}{Nm+1}$,
it follows that
$$
\dyle \iint_{\Omega_T\setminus A} |\n T_k(u_n)-\n T_k(u)|^s \psi \to 0\mbox{ as  }n\to \infty.
$$
To conclude the proof, it is sufficient to show that $ \dyle\int_A|\n T_k(u_n)|^s \psi \to 0\hbox{  as   }   n\to\infty.$

Since $T_k(u_n)\to T_k(u)\mbox{  strongly in } L^\s(\O_T)$ for all $\s>1$, then by Egorov's Lemma,  for every $\epsilon>0$, there exists a
measurable set $B_\e$  such that $|B_\e|\le \e$ an $T_k(u_n)\to T_k(u)$ uniformly in $\Omega_T\setminus B_\e$. Then
$$
\begin{array}{lll}
\dyle \iint_{\{T_k(u)=0\}}|\n T_k(u_n)|^s \psi &=& \dyle \iint_{\{T_k(u)=0\}\cap B_\e}|\n T_k(u_n)|^s+ \iint_{\{T_k(u)=0\}\cap \O\backslash
B_\e}|\n T_k(u_n)|^s \psi\\ &=& I_1+I_2
\end{array}
$$
Using the fact that $\{|\n T_k(u_n)|\}_n$ is uniformly bounded in $L^1$, and choosing $s<1$, we find that
$$
I_1\le C |B_\e|^{1-s}.
$$
By using the uniform convergence of  $T_k(u_n)$ in $\Omega_T\setminus B_\e$,  we obtain

\begin{equation}\label{gg}
I_2 \le  \dyle \iint_{\{T_k(u_n)\le M\}\cap \O\backslash B_\e}|\nabla T_k(u_n)|^s \psi \le(M+\frac 1n)^{{a s}} \iint_{\{T_k(u_n)\le M\}\cap
\O\backslash B_\e}\Big(\frac{|\n T_k(u_n)|}{(u_n+\frac 1n)^{a}}\Big)^s\psi
\end{equation}
where $a>0$. Since the estimate \eqref{UAM1} holds, it is sufficient pick up $s< 1$  and $a>0$  such that
 $$
 \iint_{\{T_k(u_n)\le M\}\cap \O\backslash
B_\e}\Big(\frac{|\n T_k(u_n)|}{(u_n+\frac 1n)^{a}}\Big)^s \psi \le C.
$$
Taking limits for  $M\to 0$ in \eqref{gg} the result follows. \cqd

\begin{remark}\label{sing}
For the case where $0<m\le \dfrac{(N-2)_+}{N}$, the  difficulty
is to show the strong convergence of the sequence $\{u_n\}_n$ in $L^1(\O_T)$.
This can be proved by assuming additional hypotheses on the structure of the measure $\mu$, see
Section \ref{sec:4}. Once proving this strong convergence, the result of the Theorem \eqref{last00} holds with the same
conclusions.
\end{remark}

\section{ The Porous Medium Equation with gradient term}
The main goal of this section is to  prove existence of solution to Problem \eqref{problemaorigen}. We start by obtaining {\it a priori
estimates} for the truncated problems, in order to be able to apply the results for an associated  {\it elliptic-parabolic} problem.

More precisely, the proof of existence follows the following steps.
\begin{enumerate}
\item  We prove some a priori estimates that allow us to show that the right hand side of the truncated problems converge weak-* to a Radon measure.
\item We transform in a natural way the problem to an {\it elliptic-parabolic} problem.
\item By using the results of Theorem \ref{last00} and compactness arguments, we identify the measure limit as the second member of the Problem
\eqref{problemaorigen}.
\end{enumerate}

We divide the section in two parts according to the values of $m$.

\subsection{ The case $1<m\le 2$}\label{sec:2}
 More precisely, consider the
problem
\begin{equation}\label{main0}
\left\{
\begin{array}{rcll}
 u_t-\D u^m&=&|\n u|^q+\, f(x,t)&\mbox{ in }\O_T\equiv \O\times
 (0,T),\\
u(x,t)&\ge &0 &\inn \O_{T},\\ u(x,t)&=&0 &\onn\p \O\times (0,T),\\ u(x,0)&=&u_0(x) &\mbox{ if }x\in\O,
\end{array}
\right.
\end{equation}
where $m>1$, $q\le 2$, $\Omega\subset \mathbb{R}^N$ a bounded domain, $f$ and $u_0$  nonnegative functions under suitable hypotheses given
below.

We will use as starting point the results in \cite{BMP1} for bounded data, $f\in L^\infty(\O_T)$ and $u_0\in L^\infty(\O)$. Since $1<m\le 2$
and $1\le q\le 2$ we will be able to obtain {\it a priori}  estimates in the framework of \cite{LKK} and \cite{KL}, where a priori estimates
are obtained to analyze the behavior of {\it viscosity supersolution}, to the porous medium equation. See  \cite{KL} and \cite{LKK} for the
details concerning to this framework.

More precisely we have the next theorem.

\begin{Theorem}\label{exis}
Assume that $1<m\le 2$ and $q\le 2$, then
\begin{enumerate}
\item If $q'(m-1)>2$, $u_0\in L^{1+\theta}(\O)$ and $f\in L^{1+\frac{2\theta}{mN}}(0,T;L^{\frac{(\theta+m)N}{mN+2\theta}}(\O))$ where
    $\theta\ge 2-m$. Then  problem \eqref{main0} has a  distributional solution.

\item If $q'(m-1)\le 2$
\begin{enumerate}
\item If  $q<m$, problem \eqref{main0} has a solution for all $f, u_0$ as in the first case. \item If $m\le q\le 2$, then problem
    \eqref{main0} has a solution if $e^{\a u_0}\in L^1(\O)$ for some $\a>0$ and $f\in L^{r}(0,T;L^{s}(\O))$ where $1<r<\infty,
    s>\frac N2$ and $\frac 1r+\frac{N}{2s}=1$.
\end{enumerate}
\end{enumerate}
\end{Theorem}

\begin{pf} \textsc{Step 1. A priori estimates.}

We prove separately each case.

{\bf (I). $q'(m-1)>2$}

Assume that $q'(m-1)>2$ and fixed $\theta >2-m$, since $m>1$, then $q<2$. Let $u_0\in L^{1+\theta}(\O)$ and $f\in
L^{1+\frac{2\theta}{mN}}(0,T;L^{\frac{(\theta+m)N}{mN+2\theta}}(\O))$, then  there exist sequences, $\{f_n\}_n, \{u_{0n}\}_n$ such that
$f_n\in L^\infty(\O_T),u_{0n}\in L^\infty(\O)$, $u_{0n}\uparrow u_0$ in $L^{1+\theta}(\O)$ and $f_n\uparrow f$ in $
L^{1+\frac{2\theta}{mN}}(0,T;L^{\frac{(\theta+m)N}{mN+2\theta}}(\O))$.

Define $u_n$, to be the bounded solution of the approximated problem
\begin{equation}\label{ap000}
\left\{
\begin{array}{rcll}
 u_{nt}-\Di(m(u_n+\frac 1n)^{m-1}\n u_n)&=&\dfrac{|\n u_{n}|^q}{|\n
 u_{n}|^q+\frac 1n}+f_n &\mbox{ in
}\O_T,\\ u_n(x,t)&=&0 &\onn\p \O\times (0,T),\\ u_n(x,0)&=&u_{0n}(x)) &\mbox{ if }x\in\O.
\end{array}
\right.
\end{equation}
Notice that the existence and the boundedness of $u_n$ follow using the result of \cite{BMP1}.

Taking $(G_k(u_n))^\theta$ as a test function in \eqref{ap000}, with $\theta>2-m$, it follows that
\begin{equation}\label{M}
\begin{array}{lll}
&\dfrac{d}{dt}\dfrac{1}{\theta +1}\dyle \io (G_k(u_n))^{\theta+1} dx+m\theta\io u_n^{m-1}(G_k(u_n))^{\theta-1}|\nabla u_n|^2dx\le
\\ &\dyle \io (G_k(u_n))^\theta|\n u_n|^q dx +\io
f_n (G_k(u_n))^\theta dx.
\end{array}
\end{equation}
Using the H\"{o}lder inequality,
\begin{equation}\label{JJ} \io (G_k(u_n))^\theta|\n
u_n|^q dx\le \epsilon \io (G_k(u_n))^{m+\theta-2}|\nabla u_n|^2dx +c(\epsilon )\io (G_k(u_n))^{\theta+\frac{q(2-m)}{2-q}}dx.
\end{equation}
Since $q'(m-1)>2$, then $\theta+\frac{q(2-m)}{2-q}< \theta +1$. Thus
$$
\io (G_k(u_n))^\theta|\n u_n|^q dx\le \epsilon  \io (G_k(u_n))^{m+\theta-2}|\nabla u_n|^2dx +c(\epsilon )\io (G_k(u_n))^{\theta+1}dx +C(\O).
$$
We deal with the last term in \eqref{M}. Using H\"older, Young and Sobolev inequalities we reach that
$$
\begin{array}{lll}
\dyle\io f_n (G_k(u_n))^\theta dx &\le &\dyle\epsilon  \Big(\io (G_k(u_n))^{\frac{(m+\theta)N}{N-2}}dx\Big)^{\frac{N-2}{N}} +C(\epsilon )
\Big(\io f^{\frac{(m+\theta)N}{mN+2\theta}}dx\Big)^{\frac{mN+2\theta}{Nm}}\\ &\le & \dyle\frac{\epsilon }{S}\io |\n
(G_k(u))^{\frac{m+\theta}{2}}|dx + C(\epsilon ) \Big(\io f^{\frac{(m+\theta)N}{mN+2\theta}}dx\Big)^{\frac{mN+2\theta}{Nm}}.
\end{array}
$$
Choosing $\epsilon $ small enough, it follows that
\begin{equation}\label{mm}
\begin{array}{lll}
&\dfrac{d}{dt}\dfrac{1}{\theta +1}\dyle \io (G_k(u_n))^{\theta+1} dx+c\io (G_k(u_n))^{\theta+m-2}|\nabla G_k(u_n)|^2dx\le
\\ &\dyle c(\epsilon )\io
(G_k(u_n))^{\theta+1}dx +C(\epsilon ) \Big(\io f^{\frac{(m+\theta)N}{mN+2\theta}}dx\Big)^{\frac{mN+2\theta}{Nm}}+ C(\O).
\end{array}
\end{equation}
Integrating in  time and using Gronwall's lemma there results that
\begin{equation}\label{m1}
\begin{array}{lll}
&C\dfrac{1}{\theta +1}\dyle \io (G_k(u_n))^{\theta+1} dx+c\iint_{\O_T}(G_k(u_n))^{\theta+m-2}|\nabla G_k(u_n)|^2dx\le
\\ &\dyle \dfrac{1}{\theta +1}\dyle \io u_0^{\theta+1} dx+
\int_0^T\Big(\io f^{\frac{(m+\theta)N}{mN+2\theta}}dx\Big)^{\frac{mN+2\theta}{Nm}}+ C(\O,T)
\end{array}
\end{equation}
Now, using $T_k(u_n)$ as a test function in the problem of $u_n$, we reach that
\begin{equation}\label{mmm}
\begin{array}{lll}
&\dyle \dyle \io \Theta_k(u_n)dx+m\iint_{\O_T}(u_n+\frac 1n)^{m-1}|\nabla T_k(u_n)|^2 \\ &\le \dyle \iint_{\O_T} T_k(u_n)|\nabla u_n|^qdxdt
+k\iint_{\O_T} f
\\
& \le \dyle \iint_{\{u_n\le \s\}}T_k(u_n)|\nabla u_n|^qdxdt +k\iint_{\{u_n\ge \s\}}|\nabla u_n|^qdxdt +k\iint_{\O_T} f
\end{array}
\end{equation}
where $\Theta_k(s)=\dyle\int_0^sT_k(\s)d\s$. From \eqref{m1} we obtain that $ \dyle \iint_{\{u_n\ge \s\}}|\nabla u_n|^qdxdt\le C(\s,f), $
thus
$$
\io \Theta_k(u_n)dx+m\iint_{\O_T}(T_k(u_n))^{m-1}|\nabla T_k(u_n)|^2 \le \iint_{\{u_n\le \s\}}T_k(u_n)|\nabla u_n|^qdxdt +C(k,a,f).
$$
Notice that, choosing $\s<<k$ small, we get the existence of $C(k)>>1$ such that $T_k^{m-1}(s)\ge C(k)T_k(s), 0\le s\le \s$, hence, using
Young's inequality (if $q<2$), there results that
\begin{equation}\label{m2}
\begin{array}{lll}
\dyle \io \Theta_k(u_n)dx+c\iint_{\O_T} T_k^{m-1}(u_n)|\nabla T_k(u_n)|^2 \le C(\O,T,k).
\end{array}
\end{equation}

Therefore, combining \eqref{m1} and \eqref{m2} we conclude that
\begin{enumerate}
\item $\{u^{1+\theta}_n\}_n$ is bounded in $L^\infty(0,T;L^1(\O))$. \item $\{(G_k(u))^{\frac{\theta+m}{2}}_n\}_n$ is bounded in
    $L^2(0,T;\sob(\O))$ and then using Poincar\'e inequality, it follows that $\{u^{m+\theta}_n\}_n$ is bounded in $L^1(\O_T)$.
\end{enumerate}

We claim that $\{|\nabla T_k(u_n)|\}_n$ is bounded in $L^2(\O_T)$ if $m<2$, while, $\{|\nabla T_k(u_n)|\}_n$ is bounded in $L^2(\O_T, \d(x))$
if $m=2$ where $\d(x)\equiv \text{dist}(x,\p\O)$.

If $m<2$, then using $w_n\equiv e^{\frac{c}{m(2-m)}(T_k(u_n)+\frac 1n)^{2-m}}-e^{\frac{c}{m(2-m)}(\frac 1n)^{2-m}}$ as a test function in
\eqref{ap000}, then we get
$$
\begin{array}{lll}
&\dyle \io L_n(u_n)dx+c\iint_{\O_T}e^{\frac{c}{m(2-m)}(T_k(u_n)+\frac 1n)^{2-m}}|\nabla T_k(u_n)|^2dx dt\le
\\ &\dyle \io L_n(u_{0})dx+\iint_{\O_T}
e^{\frac{c}{m(2-m)}(T_k(u_n)+\frac 1n)^{2-m}}|\nabla u_n|^qdx dt +C(k)\iint_{\O_T}fdx
\end{array}
$$
where $L_n(s)=\dyle\int_0^s \Big(e^{\frac{c}{m(2-m)}(T_k(s)+\frac 1n)^{2-m}}-e^{\frac{c}{m(2-m)}(\frac 1n)^{2-m}}\Big)ds$. Notice that
$L_n(s)\le C(k)s$

Thus choosing $c>>1$ and using estimate \eqref{m1}
 on $G_k(u_n)$  there results that
\begin{equation}\label{m222}
\dyle\iint_{\O_T}|\nabla T_k(u_n)|^2dx dt\le C(k).
\end{equation} Therefore the claim follows in this case.

Assume that $m=2$ and consider $\varrho$, defined in \eqref{euro000}. Fixed $0<\a<1$, to be chosen later,  using $\dfrac{\varrho}{(u_n+\frac
1n)^\a}$ as a test function in \eqref{ap000}, we get
\begin{eqnarray*}
&\dyle\frac{1}{1-\a}\io (u_n+\frac 1n)^{1-\a}\varrho dx +\frac{1}{2-\a} \iint_{\O_T}\Big((u_n+\frac 1n)^{2-\a}-(\frac 1n)^{2-\a}\Big)=\\
&\dyle \a \iint_{\O_T}\dfrac{|\n u_n|^2}{(u_n+\frac 1n)^{\a}}\varrho\, + \iint_{\O_T}\dfrac{\varrho}{(u_n+\frac 1n)^\a}\dfrac{|\n
u_{n}|^q}{|\n
 u_{n}|^q+\frac 1n} +\\
 &\dyle \iint_{\O_T}\dfrac{f_n\varrho}{(u_n+\frac
1n)^\a} +\dyle\frac{1}{1-\a}\io (u_{0n}+\frac 1n)^{1-\a}\varrho dx.
\end{eqnarray*}
Choosing $\a$ such that $1-\theta<\a<1$, then from \eqref{m1}, it follows that the first term in the above identity is uniformly bounded in
$n$. Thus
$$
\a \iint_{\O_T}\dfrac{|\n u_n|^2}{(u_n+\frac 1n)^{\a}}\varrho\,\le C
$$
and then
$$
 \iint_{\O_T}|\n T_k(u_n)|^2\varrho \le C(k).
$$
Since $\phi\backsimeq \d(x)$,  the claim follows in this case.

 Combining the result of the claim and estimates \eqref{m1} and
\eqref{m2}, we get easily that
\begin{enumerate}
\item $\{|\n u_{n}|^q\}_n$ is bounded in $L^1(\O_T)$ if $m<2$, \item $\{|\n u_{n}|^q\}_n$ is bounded in $L^1_{loc}(\O_T)$ if $m=2$.
\end{enumerate}
Notice that in booth cases we can prove that $\{u_n^{\theta+m-2}|\n u_n|^2\}_n$ is bounded in $L^1(\O_T)$.

{\bf (II) $q'(m-1)\le 2$ and $q<m$.} We deal now with the case $2)-a)$. Assume that $q'(m-1)\le 2$ and $q<m$, then the result follows using
the same kind of computations as in the first case, the main difficulty is to estimate the second term in \eqref{JJ} where we use
Poincar\'{e} inequality. This is possible using the fact that $\theta+\frac{q(2-m)}{2-q}<m+\theta.$

Since $\theta+\frac{q(2-m)}{2-q}<m+\theta$, then using H\"older and Poincar\'{e} inequalities there result that
\begin{eqnarray*}
\dyle\io (G_k(u_n))^{\theta+\frac{q(2-m)}{2-q}}dx &\le & \e \io (G_k(u_n))^{m+\theta}dx+ C(\e,\O)\\ &\le & \frac{\e}{\l_1}\io |\n
(G_k(u_n))^{\frac{m+\theta}{2}}|^2dx +C(\e,\O).
\end{eqnarray*}
Choosing $\e$ small enough and going back to estimate \eqref{m1}, it follows that
\begin{equation}\label{mm1}
\begin{array}{lll}
&C\dfrac{1}{\theta +1}\dyle \io (G_k(u_n))^{\theta+1} dx+c\iint_{\O_T}(G_k(u_n))^{\theta+m-2}|\nabla G_k(u_n)|^2dx\le C(u_0,f,T, \O).
\end{array}
\end{equation}

{\bf (III) $q'(m-1)\le 2$ and $m\le q$.}

Consider now the case where $q'(m-1)\le 2$ and $m\le q$. The existence result in this case follows using the same arguments as in
\cite{DGLS}, for the readers convenience we include here some details.

Using $e^{\a(G_k(u_n))}-1$, with $\a>0$, as a test function in the approximated problem  for $u_n$, and calling $H_k(s)=\displaystyle\int_0^s
(e^{\a(G_k(\s))}-1)d\s$, we reach that
\begin{equation}\label{HU}
\begin{array}{lll}
&\dfrac{d}{dt}\dyle \io H_k(u_n)dx+m\a\io u_n^{m-1}e^{\a(G_k(u_n))}|\nabla G_k(u_n)|^2dx =
\\ &\dyle \io (e^{\a(G_k(u_n))}-1)|\n u_n|^q dx +\io
f_n (e^{\a(G_k(u_n))}-1)dx.
\end{array}
\end{equation}
Without loss of generality we can assume that $k\ge 1$.
\begin{itemize}
\item[a)]  If $q=2$, choosing $m\a>1$, it follows that
$$
\begin{array}{lll}
& \dyle \dfrac{d}{dt}\dyle \io H_k(u_n)dx+m\a\io (u_n^{m-1}-1)e^{\a(G_k(u_n))}|\nabla u_n|^2dx +\io |\nabla G_k(u_n)|^2dx\le \\ &\dyle
\io f_n (e^{\a(G_k(u_n))}-1)dx.
\end{array}
$$
\item[b)] If $q<2$ by using Young inequality we obtain that
$$
\dyle \io (e^{\a(G_k(u_n))}-1)|\n u_n|^q dx\le \e\dyle \io e^{\a(G_k(u_n))}|\n u_n|^2dx + C(\e)\io (e^{\a(G_k(u_n))}-1)dx.
$$
Hence
$$
\begin{array}{lll}
&\dyle \dfrac{d}{dt}\dyle \io H_k(u_n)dx+m\a\io (u_n^{m-1}-\e)e^{\a(G_k(u_n))}|\nabla u_n|^2dx +\io |\nabla G_k(u_n)|^2dx\le \\ &\dyle
C(\e)\io (e^{\a(G_k(u_n))}-1)dx+\io f_n (e^{\a(G_k(u_n))}-1)dx. \end{array}
$$
Fixed $k\ge 2$, then integrating in $[0,\t]$ and taking the maximum on $\t$ it follows that
\begin{equation}\label{qmenor2}
\begin{array}{ll}
&\dyle \sup_{\t\in [0,T]}\io H_k(u_n(x,\t))dx+C\iint_{\O_T} u_n^{m-1}e^{\a(G_k(u_n))}|\nabla u_n|^2\le\\ & \dyle \io H_k(u_{0n}(x))dx+
C\dyle \iint_{\O_T}(e^{\a(G_k(u_n))}-1)+\iint_{\O_T} f_n (e^{\a(G_k(u_n))}-1)\mbox{ if }\, q<2.
\end{array}
\end{equation}
As a consequence we find that
\begin{equation}\label{laqueuso}
\iint_{\O_T} |\nabla G_k u_n|^2\le C.
\end{equation}

Moreover,
\begin{equation}\label{qigual2}
\begin{array}{ll}
&\dyle \sup_{\t\in [0,T]}\io H_k(u_n(x,\t))dx+C\iint_{\O_T} u_n^{m-1}e^{\a(G_k(u_n))}|\nabla u_n|^2\le \\ & \dyle \io
H_k(u_{0n}(x))dx+\iint_{\O_T}f_n (e^{\a(G_k(u_n))}-1)\mbox{ if }q=2.
\end{array}
\end{equation}

Let us analyze the term $\dyle\iint_{\O_T} f_n (e^{\a(G_k(u_n))}-1)$.
$$
\begin{array}{lll}
\dyle \iint_{\O_T} f_n (e^{\a(G_k(u_n))}-1)&\le &\dyle \iint_{\O_T}f (e^{\frac{\a}{2}(G_k(u_n))}-1)^2 + C \\ \\ &\le &
||f||_{r,s}||(e^{\frac{\a}{2}(G_k(u_n))}-1)||^2_{r',s'}+C.
\end{array}
$$
For simplicity of notation we set $w_n=e^{\frac{\a}{2}(G_k(u_n))}-1$. Using the Gagliardo Nirenberg inequalities there results that
$$
||w_n||^2_{r',s'}\le C||w_n||^{\frac 2r}_{\infty,2}\Big( \iint_{\O_T}|\n w_n|^2\Big)^{\frac 1r'}\le C\Big(\sup_{\t\in [0,T]}\io
w_n^2dx\Big)^{\frac 1r}\Big( \iint_{\O_T}|\n w_n|^2\Big)^{\frac 1r'}.
$$
Hence using Young's inequality we reach
$$
\begin{array}{lll}
\dyle\iint_{\O_T} f_n (e^{\a(G_k(u_n))}-1)dx &\le &\dyle C(\e)||f||^{r'}_{r,s}\iint_{\O_T}|\n w_n|^2+ \e (\sup_{\t\in [0,T]}\io
w_n^2dx\Big)+C.
\end{array}
$$
Notice that $H_k(u_n)\ge c_1 w_n^2 -c_2$, then choosing $\e$ small it follows that
$$
\begin{array}{lll}
&\dyle \sup_{\t\in [0,T]}\io H_k(u_n(x,\t))dx+c(m\a-\e) \iint_{\O_T} u_n^{m-1}|\nabla w_n|^2dx\le\\ & \dyle C(\e)||f||^{r'}_{r,s}
\iint_{\O_T} |\n w_n|^2+C(T).
\end{array}
$$
If $||f||^{r'}_{r,s}$ is sufficiently small we get
$$
\dyle \sup_{\t\in [0,T]}\io H_k(u_n(x,\t))dx+c(m\a-\e)\iint_{\O_T} u_n^{m-1}|\nabla w_n|^2\le C.
$$
If not, then we can choose $t_1<T$ such that $||f||^{r'}_{r,s}$ is sufficiently small, then
$$
\dyle \sup_{\t\in [0,t_1]}\io H_k(u_n(x,\t))dx+c(m\a-\e)\int_0^{t_1}\io u_n^{m-1}|\nabla w_n|^2dx\le C.
$$
Then the general result follows by iteration. Hence we conclude that
$$
\dyle \sup_{\t\in [0,T]}\io H_k(u_n(x,\t))dx+c\iint_{\O_T} u_n^{m-1}e^{\a(G_k(u_n))}|\nabla u_n|^2\le C(\O,T).
$$
Now, taking $T_k(u_n)$ and using the previous estimate, it follows that
\begin{equation}\label{m22}
\begin{array}{lll}
\dyle \io \Theta_k(u_n)dx+c\iint_{\O_T} T_k^{m-1}(u_n)|\nabla T_k(u_n)|^2 \le k\iint_{\O_T} f +C(\O,T,k).
\end{array}
\end{equation}
Thus there results that $|\n u_n|^q+f_n$ is bounded in $L^1(\O_T)$.
\end{itemize}

\textsc{Step 2. Passage to the limit.}

To obtain the existence of solution we need to prove  that

\begin{equation}\label{ALG0}
\dfrac{|\n u_{n}|^q}{1+\frac 1n |\n u_{n}|^q}+ f_n\to |\n u|^q +f \mbox{ strongly in } L^1_{loc}(\O_T).
\end{equation}

{\it Claim.-}  The following inequality holds,
\begin{equation}\label{sss} \int_{\{u_n\le M\}} \dfrac{|\n
u_n|^q}{(1+\frac 1n |\n u_n|)^q(u_n+\frac 1n)^{s}}\varrho\le C
\end{equation}
for all $s<1$, where  $\varrho$ is the solution of \eqref{euro000}.

To prove the claim, consider
 $\dfrac{\varrho}{(u_n+\frac 1n)^{s}}$, with $s<1$, as a test function in \eqref{ap000}. Therefore
\begin{eqnarray*}
&\dyle\iint_{\O_T} (u_n)_t (u_n+\frac 1n)^{-s}\varrho   +m\iint_{\O_T} (u_n+\frac 1n)^{m-1-s}\n u_n\n \varrho \ge \\ &\dyle s\iint_{\O_T}
\dfrac{|\n u_n|^2}{(u_n+\frac 1n)^{s+1}}\varrho +\iint_{\O_T} \dfrac{|\n u_{n}|^q}{(1+\frac 1n |\n u_{n}|^q)(u_n+\frac 1n)^{s}}\varrho
\end{eqnarray*}
and then
$$
\iint_{\O_T} \dfrac{|\n u_{n}|^q}{(1+\frac 1n |\n
 u_{n}|^q)(u_n+\frac 1n)^{s}}\varrho \le \io\frac{1}{1-s}(u_n+\frac 1n)^{1-s}\varrho dx+
\frac{1}{m-s}\iint_{\O_T} (u_n+\frac 1n)^{m-s}
$$

Since $s<1$  we obtain
\begin{equation}\label{ff}
\iint_{\O_T} \dfrac{|\n u_{n}|^q}{(1+\frac 1n |\n
 u_{n}|^q)(u_n+\frac 1n)^{s}}\varrho \le C,
\end{equation}
and the claim follows.

By using the {\it a priori} estimates of the first step, there results that
$$
\dfrac{|\n u_{n}|^q}{1+\frac 1n |\n
 u_{n}|^q}+f_n\equiv h_n \mbox{ is bounded in
 }L^1(\O_T)(\mbox{Resp. in }L^1_{loc}(\O_T)) \mbox{ if }m<2(\mbox{ Resp. if  $m=2$}).
$$
Define now
$$
v_n=(u_n+ \frac 1n)^m-(\frac 1n)^m,
$$
then $v_n$  solves the problem,
\begin{equation}\label{eq:pbapp1}
\left\{
\begin{array}{rcll}
b(v_n)_t-\Delta v_n&=&h_n,\, &(x,t)\in \O_T,\\ v_n(x,t)&=&0,\, &(x,t)\in\partial\Omega\times(0,T),\\ v_n(x,0)&=&\varphi^{-1}(T_n(u_0(x))),\,
&x\in\Omega,
\end{array}
\right.
\end{equation}
where $b(s)=(s+(\frac 1n)^m)^{\frac 1m}-\frac 1n$   satisfies the hypotheses $(B)$.

Thus by Propositions \ref{stima} and \ref{convgr}, we have that
\begin{equation}\label{ttt} \nabla v_n\to \nabla v\mbox{  strongly in } L^\sigma(\O_T)\mbox{  for all }1\le
\sigma<1+\frac 1{1+Nm}.
\end{equation}

To conclude we only need to prove \eqref{ALG0}.

The case $q=2$ is treated in \cite{DGLS} by using some kind of exponential change of variables.
We deal with the case $q<2$.

Since $\{|\n T_k(u_n)|\}_n$ is bounded in $L^2(\O_T)$  and,  by Theorem \ref{last00}, $\n u_n\to \n u$ a.e. in $\O_T$, then
$$\n T_k(u_n)\to \n T_k(u_n) \hbox{ strongly in }  L^q(\O_T).$$

We will use  Vitali's lemma. Consider $E\subset\Omega_T$,  a measurable set, then we write,
$$
\begin{array}{lll}
\dyle \iint_E |\n u_n|^q dx dt &= &\dyle \iint_{E\cap\{u_n<k\}} |\n T_k(u_n)|^q dx dt+\iint_{E\cap\{u_n\ge k\}}|\n u_n|^qdx dt.
\end{array}
$$
By using the strong convergence of the truncations, we have
$$\iint_E |\n T_k u_n|^q dx dt\to \iint_E |\n T_k u|^q dx dt.$$

We deal with the last term in the right hand side. By \eqref{m1}, \eqref{mm1} and \eqref{laqueuso} we have  in all the cases,
$$\iint_{\{u_n\ge k\}} |\n u_n|^2 \le C.$$
Then
$$\iint_{E\cap\{u_n\ge k\}}|\n u_n|^q dx dt\le C \Big(\iint_{\{u_n\ge k\}} |\n u_n|^2 \Big)^{\frac
12}|\{u_n\ge k\}|^{\frac 12}\le C |\{u_n\ge k\}|^{\frac 12}
$$
It is clear that $|\{u_n\ge k\}|\to 0$ as $k\to \infty$ uniformly in $n$. Hence the result follows using Vitali's lemma.

If $m=2$, then we can repeat the same arguments  above  to handle the term $|\n u_n|^q\psi$, for  $\psi\in \mathcal{C}^\infty_0(\O_T)$.

Therefore in   both cases we reach
$$
\dfrac{|\n u_{n}|^q}{1+\frac 1n |\n
 u_{n}|^q}\to |\n u|^q \mbox{ strongly in  }L^1_{loc}(\O_T)
$$
and  the existence result follows.
\end{pf}

\begin{remark}\label{Cont}
Notice that $1+\theta$, $1+\frac{2\theta}{mN}$, $\frac{(\theta+m)N}{mN+2\theta}\to 1$ as $\theta \to 0$. Then fixed $q<2$, for all $\e>0$,
there exists $1<m(\e)<2$ such that if $m>m(\e)$, then problem \eqref{main0} has a nonnegative solution for all $u_0\in L^{1+\e}(\O)$ and
$f\in L^{1+\e}(\O_T)$. This motivate the existence result studied in the next subsection.

\end{remark}

\subsection{ The case $m>2$: $L^1$ data}\label{sec:3}

\

In the elliptic case  if $q(\frac{1}{m}-1)<-1$,  then existence result holds for all $L^1$ data, without restriction on its size,
see \cite{ADPW}.

The goal of this subsection is to consider the case $m>2$, that implies the above condition. In particular, we can also see the next result
as a slight improvement of the result obtained  in the elliptic case.

By using suitable a priori estimates, as in the elliptic case, we will prove that problem \eqref{main0} has a distributional solution for all
$f\in L^1(\O_T)$ and $u_0\in L^1(\O)$. The main existence result is the following.

\begin{Theorem}\label{3}
Let  $f,u_0$ be such that $f\in L^1(\O_T)$ and $u_0\in L^1(\O)$. Assume $1<q\le 2$ and
 $m>2$, then problem \eqref{main0} has a
distributional solution $u$ such that $|\nabla u^m|\in L^\sigma_{loc}(\O_T)$
 for all $1\le\sigma<1+\frac{1}{Nm+1}$.
 \end{Theorem}
\begin{pf}
We will consider separately the cases $q<2$ and $q=2$.

Let begin by the case $q<2$. Define $u_n$ to be a solution to the approximated problem
\begin{equation}\label{ap01}
\left\{
\begin{array}{rclll}
 u_{nt}-\Di(m(u_n+\frac 1n)^{m-1}\n u_n)&=&
\dfrac{|\n u_{n}|^q}{1+\frac{1}{n}|\n u_{n}|^q} +T_n(f) &\mbox{ in }\O_T,\\ u_n(x,t)&=&0 &\onn\p \O\times (0,T),\\ u_n(x,0)&=&T_n(u_0(x))
&\mbox{ if }x\in\O.
\end{array}
\right.
\end{equation}
Using $e^{-\frac{c}{(m-2)(u_n+\frac 1n)^{m-2}}}\varrho$,  where $\varrho$ is defined in \eqref{euro000}, as a test function in \eqref{ap01},
it follows that
$$
\begin{array}{lll}
&\dyle \io D_n(u_n)\varrho dx+c\iint_{\O_T} e^{-\frac{c}{(m-2)(u_n+\frac 1n)^{m-2}}}|\nabla u_n|^2\varrho dx dt+\iint_{\O_T} K_n(u_n)dx
dt\le
\\ &\dyle \io D_n(u_{n0})\varrho dx+\iint_{\O_T}
e^{-\frac{c}{(m-2)(u_n+\frac 1n)^{m-2}}}|\n u_n|^q\varrho dx dt +\iint_{\O_T} e^{-\frac{c}{(m-2)(u_n+\frac 1n)^{m-2}}}\varrho T_n(f)dx dt
\end{array}
$$
with $D_n(s)=\dyle\int_0^s e^{-\frac{c}{(m-2)(t+\frac 1n)^{m-2}}}dt$ and $K_n(s)=\dyle\int_0^s m(s+\frac 1n)^{m-1}e^{-\frac{c}{(m-2)(t+\frac
1n)^{m-2}}}dt$

Since $e^{-\frac{c}{(m-2)s^{m-2}}}\le C$, for $s\ge 0$, then $c_1 s-c_2\le D_n(s)\le s$ and $K_n(s)\ge c_1 s^m-c_2$. Hence using  Young's
inequality,
$$
\begin{array}{lll}
&\dyle \io D_n(u_n)\varrho dx+c\iint_{\O_T} e^{-\frac{c}{(m-2)(u_n+\frac 1n)^{m-2}}}|\nabla u_n|^2\varrho dx dt +\iint_{\O_T} K_n(u_n)dx
dt\le
\\ &\dyle \io u_{n0}\varrho dx+\e\iint_{\O_T}
e^{-\frac{c}{(m-2)(u_n+\frac 1n)^{m-2}}}|\n u_n|^2\varrho dx dt+\\ &\dyle c\iint_{\O_T} e^{-\frac{c}{(m-2)(u_n+\frac 1n)^{m-2}}}\varrho dx dt
+\iint_{\O_T} fdx dt.
\end{array}
$$
Choosing $\e$ small il follows that
$$
\dyle \io D_n(u_n)\varrho dx+c\iint_{\O_T} e^{-\frac{c}{(m-2)(u_n+\frac 1n)^{m-2}}}|\nabla u_n|^2\varrho dx dt +\iint_{\O_T} K_n(u_n)dx dt\le
C.
$$
As a consequence, $\{u_n\}_n$ is bounded in $L^\infty(0,T;L^1_{loc}(\O))$, $\{u^m_n\}_n$ is bounded in $L^1(\O_T)$
 and $\{G_k(u_n)\}_n$ is bounded in $L^2(0,T;W^{1,2}_{loc}(\O))$.

 Using $T_k(u_n)\varrho$ as a test function in \eqref{ap01} we reach
 that the sequence $\{\Lambda_k(u_n)\}_n$ is bounded in the space
 $L^2(0,T;W^{1,2}_{loc}(\O))$ where $\dyle\Lambda_k(s)=\int_0^s
 (T_k(\s))^{\frac{m-1}{2}}d\s$.

In the same way we get easily that
\begin{equation}\label{qqq}
\iint_{\O_T}\Big(\dfrac{|\n u_{n}|^q}{1+\frac 1n |\n u_{\n}|^q} +T_n(f)\Big)\varrho \le C \mbox{ for all  }n.
\end{equation}
Fixed $s<\min\{1, m-2\}$,  and using $\dfrac{\varrho }{(u_n+\frac 1n)^s}$ as a test function in \eqref{ap01} we reach that
\begin{equation}\label{PCI}
\iint_{\O_T} (u_n+\frac 1n)^{m-2-\a}|\nabla u_n|^2\varrho \, dx dt+\iint_{\O_T} \dfrac{T_n(f)\:\varrho }{(u_n+\frac 1n)^\a}dx dt\le C \mbox{
 for all }n
 \end{equation}
and \begin{equation}\label{PCII} \iint_{\O_T}\Big(\dfrac{|\n u_{n}|^q}{1+\frac 1n |\n u_{n}|^q} \Big)\dfrac{\varrho }{(u_n+\frac 1n)^s}\le C
\mbox{
 for all }n.\end{equation}
 Thus $\{u_n^{m-2-s}|\n u_n|^2\}_n$ is bounded in
$L^1_{loc}(\O)$.

Notice that, if $2<m<3$, then there results that
$$
\iint_{\O_T}|\n T_k(u_n)|^2\varrho \le C\mbox{
 for all }n.
$$
We claim that $G_k(u_n)\to G_k(u)$ strongly in $L^2(0,T;W^{1,\a}_{loc}(\O))$ for all $\a<2$ and for all $k>0$.

By Theorem \ref{last00},  we have that $\n u_n\to \n u$ a.e in $\O_T$ in particular, $\n G_k(u_n)\to \n G_k(u)$  a.e in $\O_T$. Hence to get
the desired result we use Vitali's Lemma. Let $M>k$ and $\psi\in \mathcal{C}^\infty_0(\O_T)$, then for a measurable set $E\subset\Omega_T$ we
have
$$
\begin{array}{lll}
&\dyle  \iint_{E }|\n G_k(u_n)|^\a\psi=\\ &\dyle  \iint_{E\cap\{u_n<M\}} |\n G_k(u_n)|^\a\psi + \iint_{E\cap\{u_n\ge M\}}|\n
G_k(u_n)|^\a\psi=\\ &\dyle  \iint_{E\cap\{k\le u_n<M\}} |\n G_k(u_n)|^\a\psi + \iint_{E\cap\{u_n\ge M\}}|\n G_k(u_n)|^\a\psi.
\end{array}
$$
Since
\begin{equation}\label{forte_troncate}
      T_k (u^m_n) \to T_k (u^m)\qquad \hbox{ strongly in }
      L^2(0,T;W^{1,\a}_{loc}(\Om))\mbox{  for all }\a<2,
\end{equation}
then
$$
\dyle  \iint_{E\cap\{k<u_n<M\}} |\n G_k(u_n)|^\a\psi\le \frac \e2\mbox{  if   }\, |E|\le \d_\e.
$$
We deal now with the second term. Using \eqref{PCI} we have
$$ \iint_{E\cap\{u_n\ge M\}}
|\n u_n|^\a\psi dx\le C \Big( \iint_{\{u_n\ge M\}} u_n^{m-2-s}|\n u_n|^\a u^{-(m-2-s)}_n\psi dx dt\Big)\le\frac{C}{M^{m-2-s}}
$$
where $s>0$ is chosen such that $0<s<m-2$. Thus choosing $M$ large we reach that
$$ \iint_{E\cap\{u_n\ge M\}}|\n u_n|^\a
\psi dx\le\frac \e2.
$$
Therefore the strong convergence of $\{\n G_k(u_n)\}_n$ follows and the claim is proved.

As in the previous subsection, to get  the existence result we just have to prove that
\begin{equation}\label{ALG}
\dfrac{|\n u_{n}|^q}{1+\frac 1n |\n u_{n}|^q}\to |\n u|^q \mbox{ strongly in } L^1_{loc}(\O_T).
\end{equation}

From Theorem \ref{last00} and by the above estimate we reach $\dfrac{|\n u_{n}|^q}{1+\frac 1n |\n u_{n}|^q}\psi \to |\n u|^q \psi$ a.e. in
$\O_T$. Using \eqref{PCII}, we can prove that
\begin{equation}\label{PCIII}
\limsup\limits_{M\to 0}\dyle  \iint_{\{u_n\le M\}} \dfrac{|\n u_n|^q}{1+\frac 1n |\n u_n|^q}\psi dx=0  \mbox{ uniformly in  }  n.
\end{equation}
 Using the result of the last claim and from \eqref{forte_troncate} we obtain that
$$
\dfrac{|\n u_{n}|^q}{1+\frac 1n |\n u_{n}|^q}\to |\n u|^q \mbox{ strongly in } L^1_{loc}(\O_T\cap \{ u>0\}).
$$
If $|\{u=0\}|=0$, then \eqref{ALG} follows. Assume that $|\{u=0\}|>0$, since $|\n u|^q\in L^1_{loc}(\O_T)$, then we conclude that $|\n
u|^q=0$ on the set $\{ u=0\}$. Hence to finish the proof we have just to prove that
$$
\dfrac{|\n u_{n}|^q}{1+\frac 1n |\n u_{n}|^q}\to 0 \mbox{ strongly in } L^1_{loc}(\O_T\cap \{ u=0\}).
$$
This follows using \eqref{PCIII} and Egorov's theorem. Thus the existence result follows in this case.

We deal now with the case $q=2$. From the result of \cite{DD}, there exists  a bounded solution $u_n$ to the problem
\begin{equation}\label{ap0122}
\left\{
\begin{array}{rcllll}
 u_{nt}-\Di(m(u_n+\frac 1n)^{m-1}\n u_n)&=&|\n u_{n}|^2 +T_n(f) &\mbox{ in
}\O_T,\\ u_n(x,t)&=&0 &\onn\p \O\times (0,T),\\ u_n(x,0)&=&T_n(u_0(x)) &\mbox{ if }x\in\O.
\end{array}
\right.
\end{equation}
Using the same argument as in the case $q<2$, we find the same estimates for the sequence $\{u_n\}_n$, moveover, for all
$\psi\in\mathcal{C}^\infty_0(\O_T)$, $\psi\ge 0$,
\begin{equation}\label{PCI0}
\iint_{\O_T}|\n u_n|^2\psi <C\mbox{  for all }n
\end{equation}
and
\begin{equation}\label{PCI000}\iint_{\O_T}\dfrac{|\n u_{n}|^2\:\psi}{(u_n+\frac 1n)^s}\le C
\mbox{
 for all }n.\end{equation}
By using Theorem \ref{last00}, $\n u_n\to \n u$ e.a in $\O_T$  and then, as above,
\begin{equation}\label{PCI1}
T_k (u^m_n) \to T_k (u^m)\qquad \hbox{ strongly in } L^2(0,T;W^{1,\a}_{loc}(\Om))\mbox{  for all }\a<2.
\end{equation}
Now from \eqref{PCI0}, we conclude that
$$
u_n\to u\qquad \hbox{ strongly in } L^2(0,T;W^{1,\a}_{loc}(\Om))\mbox{  for all }\a<2.
$$
Define $w_n\equiv T_{2k}\big(u_n-T_h(u_n)+T_k(u_n)-(T_k(u))_\nu\big)_+$ where $(T_k(u))_\nu$ is defined as in \eqref{regg} and $\g(s)=\frac
1{m(2-m)}(s+\frac 1n)^{2-m}$. Let $h>2 k>0$ to be chosen later. It is clear that $\n w_n\equiv 0$ for $u_n>M\equiv 4k+h$.

Using $w_ne^{\g(u_n)}\psi$ as a test function in \eqref{ap0122}, it follows that
\begin{eqnarray*}
&\dyle \int_0^T \langle (u_n)_t, e^{\g(u_n)}w_n\psi \rangle dt +m\iint_{\O_T} e^{\g(u_n)}(u_n+\frac 1n)^{m-1}\psi\n T_M(u_n)\n w_n\\&\dyle
+m\iint_{\O_T} e^{\g(u_n)}(u_n+\frac 1n)^{m-1}w_n\n T_M(u_n)\n \psi\le \iint_{\O_T} T_n(f) e^{\g(u_n)}w_n\psi
\end{eqnarray*}
It is clear that
$$
\iint_{\O_T} T_n(f)e^{\g(u_n)}w_n \psi \le \o(n,\nu)
$$
and
$$
|\iint_{\O_T} e^{\g(u_n)}(u_n+\frac 1n)^{m-1}w_n\n T_M(u_n)\n \psi|\le \o(n,\nu).
$$
It is well known that
$$
\dyle \int_0^T \langle (u_n)_t, e^{\g(u_n)}w_n\psi \rangle dt\ge \o(n)+\o(\nu),
$$
see, for instance,  \cite{DGLS}.

Let analyze the term   $m\dyle\iint_{\O_T} e^{\g(u_n)}(u_n+\frac 1n)^{m-1}\n T_M(u_n)\n w_n\psi$.
$$
\begin{array}{lll}
&\dyle \iint_{\O_T}e^{\g(u_n)}(u_n+\frac 1n)^{m-1} \n T_M(u_n)\n w_n\psi \\ &= \dyle \iint_{\{u_n\le k\}} e^{\g(u_n)}(u_n+\frac 1n)^{m-1} \n
T_k(u_n)\n w_n\psi +\int_{\{u_n>k\}} e^{\g(u_n)}\psi (u_n+\frac 1n)^{m-1} \n T_M(u_n)\n w_n\\ &\ge \dyle \iint_{\O_T} e^{\g(u_n)}(u_n+\frac
1n)^{m-1} \n T_k(u_n)\n (T_k(u_n)-(T_k(u))_\nu)_+\psi \\ &\dyle- \iint_{\{u_n>k\}} e^{\g(u_n)}(u_n+\frac 1n)^{m-1} |\n T_M(u_n)||\n
(T_k(u))_\nu|\psi.
\end{array}
$$
Define
$$
\Gamma_n(s)= \left\{
\begin{array}{lll}
e^{\g(s)}(s+\frac 1n)^{m-1}&\mbox{  if  }& s\le k,\\ 0&\mbox{  if   }& s\ge k,
\end{array}
\right.
$$
then
$$
\begin{array}{lll}
&\dyle \iint_{\O_T} e^{\g(u_n)}(u_n+\frac 1n)^{m-1} \n T_k(u_n)\n (T_k(u_n)-(T_k(u))_\nu)_+\psi\\&=\dyle  \iint_{\O_T} \Gamma_n(u_n)\n
T_k(u_n)\n (T_k(u_n)-(T_k(u))_\nu)_+\psi\\ =&\dyle \iint_{\O_T} \Gamma_n(u_n)|\n (T_k(u_n)-(T_k(u))_\nu)_+|^2\psi +\o(n,\nu).
\end{array}
$$
On the other hand, for $M$ fixed,
$$
\iint_{\{u_n>k\}} e^{\g(u_n)}(u_n+\frac 1n)^{m-1} |\n T_M(u_n)||\n (T_k(u))_\nu\psi =\o(n, \nu).
$$
Hence
$$
\begin{array}{ll}
&\dyle \iint_{\O_T} \Gamma_n(u_n)|\n (T_k(u_n)-(T_k(u))_\nu)_+|^2\psi\le\dyle
 \iint_{\O_T}e^{\g(u_n)}(u_n+\frac 1n)^{m-1} \n T_M(u_n)\n w_n \psi +\\
&\dyle\iint_{\{u_n>k\}} e^{\g(u_n)}(u_n+\frac 1n)^{m-1} |\n T_M(u_n)||\n (T_k(u))_\nu \psi + \o(n,\nu).
\end{array}
$$
Thus
$$
\dyle\iint_{\O_T} \Gamma_n(u_n)|\n (T_k(u_n)-(T_k(u))_\nu)_+|^2\psi \le \o(n,\nu).
$$
In the same way we reach that
$$
\dyle\iint_{\O_T}\Gamma_n(u_n)|\n (T_k(u_n)-(T_k(u))_\nu)_-|^2\psi \le \o(n,\nu).
$$
Since $\Gamma_n(s)\ge C_1$ if $s>C_2>0$, uniformly in $n$, then
\begin{equation}\label{ZAA0}
|\n u_n|\chi_{\{ c_1<u_n<c_2\}}\to |\n u|\chi_{\{ c_1<u<c_2\}}\mbox{  strongly   in }L^2_{loc}(\O_T).
\end{equation}
We claim that
\begin{equation}\label{ZAA}
|\n G_k(u_n)|\to |\n G_k(u)| \mbox{  strongly   in }L^2_{loc}(\O_T).
\end{equation}
We again  use Vitali's lemma.  Let $\psi\in \mathcal{C}^\infty_0(\O_T)$ be such that $\psi\ge 0$. Consider $E\subset \O_T$, a
measurable set and write,
$$
\begin{array}{lll}
\dyle  \iint_{E} |\n G_k(u_n)|^2\psi &= & \dyle \iint_{E\cap \{u_n\le M\}} |\n G_k(u_n)|^2\psi+  \iint_{E\cap\{ u_n\ge M\}} |\n
G_k(u_n)|^2\psi\\ &=& \dyle  \iint_{E\cap \{k\le u_n\le M\}} |\n G_k(u_n)|^2\psi+
 \iint_{E\cap \{u_n\ge M\}} |\n G_k(u_n)|^2\psi.
\end{array}
$$
From \eqref{ZAA0}, given $\epsilon>0$ there exists $\delta>0$ such that
$$
\limsup\limits_{n\to \infty} \iint_{E\cap\{ k\le u_n\le M\}} |\n G_k(u_n)|^2\psi\le \frac \e2 \mbox{  if   }|E|\le \d.
$$
We deal now with the term $\int_{E\cap\{ u_n\ge M\}} |\n G_k(u_n)|^2\psi.$

Using \eqref{PCI}, we get, for some $\a<m-2$,
$$
 \iint_{E\cap \{u_n\ge M\}} |\n G_k(u_n)|^2\psi= \iint_{E\cap \{u_n\ge M\}}
\frac{(u_n+\frac 1n)^{m-2-\a}}{(u_n+\frac 1n)^{m-2-\a}}|\nabla u_n|^2\psi\le \frac{C}{(M+\frac 1n)^{m-2-\a}}.
$$
Thus
$$
\limsup \limits_{n\to \infty}\dyle  \iint_{E} |\n G_k(u_n)|^2\psi\le \frac \e2+\frac{C}{(M+\frac 1n)^{m-2-\a}}.
$$
Letting $M\to \infty$, we reach the claim.

In the same way and by using estimates \eqref{PCI000}, we can prove that
$$
\iint_{\O_T}|\n (T_k(u_n)-T_k(u))|^2\psi\to 0\mbox{ as  }n\to \infty.
$$
Thus
$$
|\n u_n|^2\to |\n u|^2\mbox{  strongly  in }L^1_{loc}(\O_T),
$$
and then the existence result follows.
\end{pf}

\begin{remark}
If $m=q=2$ and $f=0$,  we  can prove that the solution to problem \eqref{main0} has the {\it finite speed propagation  property}. This
follows by setting $w=\frac{2}{3}\Big(\frac{4}{5}\Big)^{\frac 25}u^{\frac 52}$, then $w$ solves
\begin{equation}\label{extinn}
     \left\{\begin{array}{rcll}
 w_t- \frac{4}{5}\Big(\frac{3}{2}\Big)^{\frac 53}
\Delta\,w^{\frac 53}&=&0 \quad&\hbox{ in } \ \Omega_T\,,\\
 w(x,t)&=&0 \quad&\hbox{ on  }
  \p \O\times (0,T)\,,\\
  w(x,0)&=&\frac{2}{3}\Big(\frac{4}{5}\Big)^{\frac 25}u^{\frac 52}_0(x)\quad&\hbox{ in  }
  \O.
 \end{array}\right.
\end{equation}

If $u_0\in L^\infty(\O)$ has a compact support, by using a convenient Barenblatt self-similar super-solution (see \cite{V}, for instance) we
obtain the {\it finite speed of propagation  property}. The inverse change of variable allow us to conclude the same result for problem
\begin{equation}\label{extinno}
     \left\{\begin{array}{rcll}
 u_t- \Delta\,u^2&=&|\nabla u|^2 \quad&\hbox{ in } \ \Omega_T\,,\\
 u(x,t)&=&0 \quad&\hbox{ on  }   \p \O\times (0,T)\,,\\
  u(x,0)&=&u_0(x)\quad&\hbox{ in  } \O \,.
 \end{array}\right.
\end{equation}
\end{remark}

\section{ The fast diffusion equation}\label{sec:4}
In this section we consider the case  $0<m<1$, usually called {\it fast diffusion equation} in the literature. We will prove the following
existence result.

\begin{Theorem}\label{fast}
Assume that $0<m<1$,  $q\le 2$ and
\begin{enumerate}
\item $f\in L^{r}(0,T;L^{s}(\O))$ where $1<r<\infty, s>\frac N2$ with $\frac1r+\frac{N}{2s}=1$ \item $e^{\a u^{2-m}_0}\in L^1(\O)$
    where,
\begin{enumerate}
\item either $\a>0$ is any positive constant if $q<2$ \item or $\a m(2-m)>2$ if $q=2$.
\end{enumerate}
\end{enumerate}
Then  problem \eqref{main0} has a distributional solution
\end{Theorem}

\begin{pf}
Let $\{f_n\}_n, \{u_{0n}\}_n$ be sequences of bounded nonnegative functions such that $u_{0n}\uparrow u_0$ and $f_n\uparrow f$.

Let $u_n$ be the bounded solution of
\begin{equation}\label{ap001}
\left\{
\begin{array}{rclll}
 u_{nt}-\Di(m(u_n+\frac 1n)^{m-1}\n u_n)&=&\dfrac{|\n u_{n}|^q}{1+\frac 1n |\n u_{n}|^q} +f_n &\mbox{ in
}\O_T,\\ u_n(x,t)&=&0 &\onn\p \O\times (0,T),\\ u_n(x,0)&=&u_{0n}(x)) &\mbox{ if }x\in\O,
\end{array}
\right.
\end{equation}
with data $(f_n,u_{0n})$. Notice that the existence and the boundedness of $u_n$ follow using the result of \cite{BMP1}.

Taking $e^{\a u^{2-m}_n}-1$, $\a>0$, as a test function in \eqref{ap001}, we find that
\begin{equation}\label{HUU}
\begin{array}{lll}
&\dfrac{d}{dt}\dyle \io H(u_n)dx+m\a(2-m)\io e^{\a u^{2-m}_n}|\nabla u_n|^2dx\le
\\ &\dyle \io (e^{\a u^{2-m}_n}-1)|\n u_n|^q dx +\io
f_n (e^{\a u^{2-m}_n}-1)dx
\end{array}
\end{equation}
where $H(s)=\dyle\int_0^s (e^{\a \s^{2-m}}-1)d\s$.

Using Young inequality, integrating in $[0,\t]$ and taking the maximum on $\t$
$$
\begin{array}{ll}
&\dyle \sup_{\t\in [0,T]}\io H(u_n(x,\t))dx+(\a m(2-m)-\e)\iint_{\O_T} e^{\a u^{2-m}_n}|\nabla u_n|^2\le\\ & C\dyle \iint_{\O_T}(e^{\a
u^{2-m}_n}-1)+\iint_{\O_T} f_n (e^{\a u^{2-m}_n}-1)dx +\io H(u_0(x))dx \mbox{ if }\,\,q<2
\end{array}
$$
and
$$
\begin{array}{ll}
&\dyle \sup_{\t\in [0,T]}\io H(u_n(x,\t))dx+(\a m(2-m)-1)\iint_{\O_T} e^{\a u^{2-m}_n}|\nabla u_n|^2dx\le \\ & C\dyle \iint_{\O_T} (e^{\a
u^{2-m}_n}-1)+\io f_n (e^{\a u^{2-m}_n}-1)dx +\io H(u_0(x))dx \mbox{ if } \,\,q=2.
\end{array}
$$
Let us analyze the last term in \eqref{HUU}.
$$
\begin{array}{lll}
\dyle\iint_{\O_T} f_n (e^{\a u^{2-m}_n}-1)dx &\le &\dyle \iint_{\O_T} f (e^{\frac{\a}{2}u^{2-m}_n}-1)^2 + C \\ &\le &
||f||_{r,s}||(e^{\frac{\a}{2}u^{2-m}_n}-1)||^2_{r',s'} +C
\end{array}
$$
We set $w_n=e^{\frac{\a}{2}u^{2-m}_n}-1$, then using the Gagliardo-Nirenberg inequality we obtain that
$$
||w_n||^2_{r',s'}\le C||w_n||^{\frac 2r}_{\infty,2}\Big(\int_0^\t\io |\n w_n|^2\Big)^{\frac 1r'}\le C\Big(\sup_{\t\in [0,T]}\io
w_n^2dx\Big)^{\frac 1r}\Big(\int_0^\t\io |\n w_n|^2\Big)^{\frac 1r'}.
$$
Thus
$$
\begin{array}{lll}
\dyle\iint_{\O_T} f_n (e^{\a u^{2-m}_n}-1)dx &\le &\dyle C(\e)||f||^{r'}_{r,s} \int_0^\t\io |\n w_n|^2+ \e (\sup_{\t\in [0,T]}\io
w_n^2dx\Big)+C.
\end{array}
$$
Using the fact that $H(u_n)\ge c_1 w_n^2 -c_2$, then choosing $\e$ small it follows that
$$
\begin{array}{ll}
&\dyle \sup_{\t\in [0,T]}\io H(u_n(x,\t))dx+c(\a m(2-m)-\e)\iint_{\O_T} |\nabla w_n|^2dx\le\\ & \dyle C(\e)||f||^{r'}_{r,s} \int_0^\t\io |\n
w_n|^2+C(\O,T).
\end{array}
$$
If $||f||^{r'}_{r,s}$ is sufficiently small we get
$$
\dyle \sup_{\t\in [0,T]}\io H(u_n(x,\t))dx+c(\a m(2-m)\e)\iint_{\O_T} |\nabla w_n|^2dx\le C.
$$
If not, then we can choose $t_1<T$ such that $||f||^{r'}_{r,s}$ is sufficiently small, then
$$
\dyle \sup_{\t\in [0,t_1]}\io H(u_n(x,\t))dx+c(\a m(2-m)-\e)\int_0^{t_1}\io |\nabla w_n|^2dx\le C.
$$
Then the general result follows by iteration. Hence we conclude that
\begin{equation}\label{hola}
\dyle \sup_{\t\in [0,T]}\io H(u_n(x,\t))dx+c\iint_{\O_T} e^{\a u^{2-m}_n}|\nabla u_n|^2 dx\le C(\O,T).
\end{equation}
Therefore we conclude that $|\n u_n|^q+f_n$ is bounded in $L^1(\O_T)$.

It is clear form the above estimate that $\{u_n\}_n$ is bounded in the spaces $L^2(0,T;\sob(\O))\cap L^\infty(0,T;L^2(\O))$ and
$\{u_{nt}\}_n$ is bounded in the spaces $L^2(0,T;W^{-1,2}(\O))+L^1(\O_T)$.
 Then,  there exists a measurable function $u\in L^2(0,T;\sob(\O))\cap L^\infty(0,T;L^2(\O))$ such that, up to a
subsequence, $u_n\to u$ strongly in $\mathcal{C}(0,T;L^2(\O))$.
Hence using by the result of Theorem \ref{last00} and Remark
\ref{sing}, we get $\n u_n\to \n u$ a.e in $\O_T$.
Thus if $q<2$, then by the previous estimate we
obtain that $|\n u_n|^q\to |\n u|^q$ strongly in $L^1(\O_T)$ and then the result follows.

We deal now with the case $q=2$. Consider $\psi\in\mathcal{C}^\infty_0(\O_T)$,  $\psi\ge 0$. By using $\dfrac{\psi}{(u_n+\frac 1n)^\d}$ as a
test function in \eqref{ap001}, where $\d<m$, there results that

\begin{equation}\label{TTT}
\iint_{\O_T} \dfrac{|\n u_n|^2}{(u_n+\frac 1n)^{2+\d-m}}\psi\le C\mbox{ for all  }n\mbox{ and for all  } n.
\end{equation}

Notice that from \eqref{TTT}, it follows that
\begin{equation}\label{DZA}
\iint_{\{u_n\le M\}} |\n u_n|^2\psi\le C(M+\frac 1n)^{2+\d-m}\mbox{ uniformely in }n,
\end{equation}
 and
\begin{equation}\label{DZA1}
\iint_{\O_T} \dfrac{|\n T_k(u_n)|^2}{(u_n+\frac 1n)^{2+\d-m}}\psi\le C\mbox{ for all  }n.
\end{equation}
Let us prove now that \begin{equation}\label{ssss} |\n T_k(u_n)|\to |\n T_k(u)|\mbox{  strongly in }L^2_{loc}(\O_T).
\end{equation}

Let $\varphi$ be a real differentiable function such that $\varphi(0)=0$ and $(2k)^{m-1}\varphi'-|\varphi|\ge C>0$. Consider $w_n\equiv
e^{-\g(T_k(u_n))}\varphi(T_n(u_n)-(T_k(u))_\nu)_+$ where $(T_k(u))_\nu$ is defined as in \eqref{regg} and
$$\g(s)=\frac 1{m(2-m)}\Big((s+\frac 1n)^{2-m}-(\frac 1n)^{2-m}\Big).$$

Using $w_ne^{\g(u_n)}\psi $ as a test function in \eqref{ap001}, it follows that
\begin{eqnarray*}
&\dyle \int_0^T \langle (u_n)_t, e^{\g(u_n)}w_n\psi \rangle dt +m\iint_{\O_T} e^{\g(u_n)}(u_n+\frac 1n)^{m-1}\n u_n\n w_n \psi\\ &\dyle
+m\iint_{\O_T} e^{\g(u_n)}(u_n+\frac 1n)^{m-1}w_n \n u_n\n \psi\dyle \le \iint_{\O_T} f e^{\g(u_n)}w_n \psi
\end{eqnarray*}
Notice that
$$
\dyle \int_0^T \langle (u_n)_t, e^{\g(u_n)}w_n\psi \rangle dt\ge \o(n)+\o(\nu).
$$
(see for instance \cite{DGLS}). Using the  hypothesis on $f$ and by \eqref{hola} it follows that
$$
\iint_{\O_T} f e^{\g(u_n)}w_n \psi \le \o(n,\nu).
$$
On the other hand,
\begin{eqnarray*}
&\dyle |\iint_{\O_T} e^{\g(u_n)}(u_n+\frac 1n)^{m-1}w_n \n u_n\n \psi|\le \\ &\dyle \Big(\iint_{\text{Supp}(\psi)} e^{\g(u_n)}(u_n+\frac
1n)^{2(m-1)}|\n u_n|^2\Big)^{\frac 12} \Big(\iint_{\text{Supp}(\psi)} e^{\g(u_n)}w_n^2|\n \psi|^2\Big)^{\frac 12}\le \\ &\dyle C
\Big(\iint_{\text{Supp}(\psi)} e^{\g(u_n)}w_n^2|\n \psi|^2\Big)^{\frac 12}.
\end{eqnarray*}
Notice that
$$ e^{\g(u_n)}w_n^2|\n \psi|^2\to 0\:\:\hbox{  a.e. in } \Omega_T \:\:\hbox{  and   } e^{\g(u_n)}w_n^2|\n \psi|^2\le C e^{\g(u_n)},$$
therefore by the Sobolev inequality and the estimates \eqref{hola} and \eqref{TTT}, we reach that
$$\Big(\iint_{\text{Supp}(\psi)} e^{\g(u_n)}w_n^2|\n
\psi|^2\Big)^{\frac 12}=\omega(n).$$

Therefore we conclude that
$$
m\iint_{\O_T} e^{\g(u_n)}(u_n+\frac 1n)^{m-1}\n u_n\n w_n \psi\le \o(n)+\o(\nu).
$$
Notice that,
\begin{eqnarray*}
&\dyle m\iint_{\O_T} e^{\g(u_n)}(u_n+\frac 1n)^{m-1}\n u_n\n w_n \psi=\\ &\dyle m\iint_{\O_T} e^{\g((u_n))-\g(T_k(u_n))}(u_n+\frac
1n)^{m-1}\n u_n\n (T_k(u_n)-T_k(u))_+\varphi'((T_k(u_n)-T_k(u))_+) \psi\\ &-\dyle m\iint_{\{u_n\le k\}} |\n
u_n|^2\varphi((T_k(u_n)-T_k(u))_+) \psi
\end{eqnarray*}
Let us analyze each term in the previous identity.

\begin{eqnarray*}
&\dyle m\iint_{\O_T} e^{\g((u_n))-\g(T_k(u_n))}(u_n+\frac 1n)^{m-1}\n u_n\n (T_k(u_n)-T_k(u))_+\varphi'((T_k(u_n)-T_k(u))_+) \psi=\\ &\dyle
\dyle m\iint_{\{u_n\le k\}} (u_n+\frac 1n)^{m-1}|\n (T_k(u_n)-(T_k(u))_\nu)_+|^2\varphi'((T_k(u_n)-T_k(u))_+) \psi\\ &\dyle  +  \dyle
m\iint_{\{u_n\le k\}} (u_n+\frac 1n)^{m-1} \n((T_k(u))_\nu))\n (T_k(u_n)-(T_k(u))_\nu)_+\varphi'((T_k(u_n)-T_k(u))_+) \psi\\ &\dyle -\dyle
m\iint_{\{u_n\ge k\}} (u_n+\frac 1n)^{m-1} \n((T_k(u))_\nu))\n (T_k(u_n)-(T_k(u))_\nu)_+\varphi'((T_k(u_n)-T_k(u))_+) \psi
\end{eqnarray*}
It is clear that
$$
|\dyle\iint_{\{u_n\ge k\}} (u_n+\frac 1n)^{m-1} \n((T_k(u))_\nu))\n (T_k(u_n)-(T_k(u))_\nu)_+\varphi'((T_k(u_n)-T_k(u))_+) \psi|\le
\o(n,\nu).$$ Now,

\begin{eqnarray*}
&\dyle\iint_{\{u_n\le k\}} (u_n+\frac 1n)^{m-1} \n((T_k(u))_\nu))\n (T_k(u_n)-(T_k(u))_\nu)_+\varphi'((T_k(u_n)-T_k(u))_+) \psi=\\ &
\dyle\iint\limits_{\{u_n\le k\}\cap \{T_k(u_n)\ge (T_k(u))_\nu\}} \psi(u_n+\frac 1n)^{m-1} \n((T_k(u))_\nu))\n
(T_k(u_n)-(T_k(u))_\nu)_+\varphi'((T_k(u_n)-T_k(u))_+)
\end{eqnarray*}
Since $m<1$, by using the dominated convergence theorem and by estimate \eqref{DZA1},
$$
\psi(u_n+\frac 1n)^{m-1} |\n((T_k(u))_\nu))|\chi_{\{T_k(u_n)\ge (T_k(u))_\nu\}} \to  u^{m-1} |\n T_k(u)|\psi \:\mbox{  strongly in
 } L^2(\O_T)\mbox{  as }n,\nu\to \infty.
$$
Thus using a duality argument, we reach
$$
\dyle m\iint_{\{u_n\le k\}} (u_n+\frac 1n)^{m-1} \n((T_k(u))_\nu))\n (T_k(u_n)-(T_k(u))_\nu)_+\varphi'((T_k(u_n)-T_k(u))_+) \psi=\o(n,\nu).
$$
It is clear that
\begin{eqnarray*}
&\dyle -m\iint_{\{u_n\le k\}} |\n u_n|^2\varphi((T_k(u_n)-T_k(u))_+) \psi\ge\\ &\dyle  -m\iint_{\{u_n\le k\}} |\n
(T_k(u_n)-(T_k(u))_\nu)_+|^2\varphi((T_k(u_n)-T_k(u))_+) \psi+\o(n)+\o(\nu). \end{eqnarray*}
 Thus combining the above estimates, there
results that
$$
m\iint_{\{u_n\le k\}}\Big((u_n+\frac 1n)^{m-1}\varphi'-\varphi)|\n (T_k(u_n)-((T_k(u))_\nu)_+)|^2\le \o(n)+\o(\nu)
$$
and by the properties of $\varphi$,
$$
\iint_{\O_T}|\n (T_k(u_n)-T_k(u))_+)|^2\psi\to 0\mbox{ as  }n\to \infty.
$$
In a similar way, we obtain
$$
\iint_{\O_T}|\n (T_k(u_n)-T_k(u))_-)|^2\psi\to 0  \mbox{ as  }n\to \infty.
$$
Hence
$$
|\n T_k(u_n)|\to |\n T_k(u)|\mbox{  strongly  in }L^2_{loc}(\O_T).
$$
Thus using estimate \eqref{TTT}, \eqref{ssss} and  Vitali's lemma, we can prove that
$$
|\n u_n|\to |\n u|\mbox{  strongly  in }L^2_{loc}(\O_T).
$$
and then the existence result follows.
\end{pf}

\subsection{Finite time extinction}\label{sub:41}
Assume that $f\equiv 0$ and $q=2$, that is, we consider the problem
\begin{equation}\label{exte}
\left\{
\begin{array}{rclll}
 u_{t}-\D u^m&=&|\n u|^2 &\mbox{ in
}\O_T,\\ u(x,t)&=&0 &\onn\p \O\times (0,T),\\ u(x,0)&=&u_0(x) &\mbox{ if }x\in\O,
\end{array}
\right.
\end{equation}
where $0<m<1$.

We will prove that the {\it regular} solutions of \eqref{exte} become zero in finite time, provided the initial datum $u_{0}\in
L^{1+\theta_0}(\O)$ for some $\theta_0>0$.

We begin by defining the meaning of {\it regular} solutions to \eqref{exte}.

\begin{Definition}\label{ddd1}
Let $ H(s)= \dyle\int_0^s e^{\frac{t^{\frac{2-m}{m}}}{m(2-m)}}dt $ and define
\begin{equation}\label{eq:fi} \dyle\beta(s)=\frac{1}{m}\int_0^s(H^{-1}(\s))^{\frac{1}{m}-1}d\s,
\end{equation}
we say that $u$ is a regular solution to problem \eqref{exte} in $\O_T$ if
$$
v\equiv H(u^m)\in
    L^2((0,T);W_0^{1,2}(\Omega))\cap \mathcal{C}([0,T]; L^2(\Omega)), \,\b(v)_t\in
    L^{2}((0,T);W^{-1,2}(\Omega))
$$
and for all $\phi\in
    L^2((0,T);W^{1,2}_0(\Omega))$ we have
    \begin{equation}\label{norma00}
    \int_0^T\langle (\b(v))_t,\phi\rangle+
    \int_0^T\int_\Omega\nabla v\cdot\nabla \phi=0.
 \end{equation}
\end{Definition}
It is clear that the existence of a regular solution follows using Theorem \ref{fast} for $q=2$, $f=0$ and the regularity of $u_0$.

We are able now to state the next result.
\begin{Theorem} \label{extinct1} Assume that $0<m<1$.
If $u$ is the regular solution of problem \eqref{exte}  in the sense of Definition \ref{ddd1}, then there exists a positive, finite time
$t_{0}$, depending on $N$,  and $u_0$ such that $u(x,t)\equiv 0$ for $t>t_{0}$.
\end{Theorem}

\begin{pf}
To get the desired result we have just to show that $v(x,t)\equiv 0$ for $t>t_{0}$. It is clear that $v$ solves
\begin{equation}\label{prob111}
     \left\{\begin{array}{rcll}
 \big(\b(v)\big)_t-\Delta\,v&=&0\quad
 &\hbox{ in } \ \Omega\times(0,T)\,,\\[5pt]
 v(x,t)&=&0 &\hbox{ on  }
  \p \O\times (0,T)\,,\\[5pt]
  v(x,0)&=&v_0(x)&\hbox{ in  } \O \,,
 \end{array}\right.
\end{equation}
with $v_0\in L^2(\O)$. Using $v^\theta$, where $\theta>0$ to be chosen later,  as a test function in \eqref{prob111}, there result that
\begin{equation}\label{z1}
\frac{d\ }{dt}\int_{\Om}\Psi(v(x,t))\,dx + c(\theta) \int_{\Om}\big|\nabla v^{\frac{\theta+1}{2}}|\,dx =0\,,
\end{equation}
where
$$
\Psi(s) = \int_{0}^{s} s^\theta (H^{-1}(\s))^{\frac{1}{m}-1}\,d\s\,.
$$
Since $$ \lim\limits_{s\to \infty}\dfrac{H^{-1}(s)}{s^\e}=0\mbox{ for all }\e>0$$ it follows that
$$
\Psi(s) \leq c(\e)s^{\theta+1+\e(\frac 1m -1)}\qquad \text{for every $s\geq0$.}
$$
Fixed $\theta$ such that $\theta+1+\e(\frac 1m -1)=1+\theta_0$, then using Sobolev's and H\"older's inequalities,
\begin{eqnarray*}
\dyle \int_{\Om}\big|\nabla v^{\frac{\theta+1}{2}}\big|^2\,dx \geq
   c_{1}(N, \theta)\bigg[\int_{\Om}\big(v^{\frac{\theta+1}{2}}\big)^{2^{*}}\,dx\bigg]^{2/2^{*}}\geq
   c_{2}(N,\theta,|\Om|)
   \bigg[\int_{\Om}(v^{a(\theta+1)})\,dx\bigg]^{1/a}
\end{eqnarray*}
where $1<a<\frac{2^*}{2}$ is chosen such that $a(\theta +1)=\theta+1+\e(\frac 1m -1)$. Hence it follows that
$$
\int_{\Om}\big|\nabla v^{\frac{\theta+1}{2}}\big|^2\,dx\geq
   c(N,\theta,|\Om|)
   \bigg[\int_{\Om}\Psi(v(x,t))\,dx\bigg]^{1/a}.
$$
Define
$$
   \xi(t)= \int_{\Om}\Psi(v(x,t))\,dx\,,
$$
then
$$
\frac{\xi'(t)}{\xi(t)^{1/a}}\leq -c_{4}<0\,.
$$
Note that by the assumption on $v_{0}$ we reach that $\xi(0)<\infty$. Integrating in $t$, one obtains
$$
   \frac{a}{a-1}\big(\xi(t)^{\frac{a-1}{a}}- \xi(0)^{\frac{a-1}{a}}\big)
   \leq
   -c_{4}t\,.
$$
Thus, as long as $\xi(t)>0$, one has
$$
   \xi(t)^{\frac{a-1}{a}}\leq \xi(0)^{\frac{a-1}{a}} - c_{4}\frac{a-1}{a}\,t\,.
$$
Therefore, $\xi(t)\equiv 0$ for $t$ large enough.
\end{pf}

\end{document}